\input ./style/arxiv-general.cfg
\documentclass[MSNbibl,number,citesort,seceqn,dvips]{arxbj}
\makeatletter
   \@ifpackageloaded{graphicx}{}{\usepackage{graphicx}}
\makeatother
\usepackage{mathbh,upgreek}

\aid{0}
\volume{22}
\issue{1}
\pubyear{2016}
\firstpage{213}
\lastpage{241}
\doi{10.3150/14-BEJ656} 
\docsubty{FLA}

\makeatletter

\newtheorem{theorem}{Theorem}[section]
\newtheorem{proposition}[theorem]{Proposition}
\newtheorem{lemma}[theorem]{Lemma}
\newremark{Remark}{Remark}
\newremark{Example}{Example}

\newcommand{\indic}{\mathbh{1}}
\newcommand{\reals}{\mathbb{R}}
\newcommand{\prob}{\mathrm{P}}
\newcommand{\expec}{\mathrm{E}}
\newcommand{\mse}{\operatorname{MSE}}
\newcommand{\argmax}{\operatorname{\arg\max}}
\newcommand{\matlab}{\textsc{Matlab}\ }

\newcommand{\Binomial}{\operatorname{Binomial}}
\newcommand{\Bernoulli}{\operatorname{Bernoulli}}
\newcommand{\eqref}[1]{(\ref{#1})}
\newcommand{\Hypergeo}{\operatorname{Hypergeo}}
\newcommand{\area}{\operatorname{area}}
\newcommand{\CFG}{\operatorname{CFG}}
\newcommand{\CFGopt}{\operatorname{CFGopt}}

\makeatother

\begin{document}
\begin{frontmatter}

\title{Polynomial Pickands functions}
\runtitle{Polynomial Pickands functions}

\begin{aug}
\author[1]{\inits{S.}\fnms{Simon}~\snm{Guillotte}\corref{}\thanksref{1}\ead[label=e1]{guillotte.simon@uqam.ca}}
\and
\author[2]{\inits{F.}\fnms{Fran\c{c}ois}~\snm{Perron}\thanksref{2}\ead[label=e2]{perronf@dms.umontreal.ca}}
\address[1]{D\'{e}partement de Math\'{e}matiques, Universit\'e du Qu\'
{e}bec \`{a} Montr\'{e}al, 201 Ave du
Pr\'{e}sident-Kennedy, Montr\'{e}al, Ou\'{e}bec H2X 3Y7, Canada. \printead{e1}}
\address[2]{D\'{e}partement de Math\'{e}matiques et Statistique, Universit\'e de Montr\'eal,
C.P. 6128, Succursale Centreville,  Montr\'{e}al, Ou\'{e}bec H3C 3J7, Canada. \printead{e2}}
\end{aug}

\received{\smonth{5} \syear{2013}}
\revised{\smonth{6} \syear{2014}}

%
\begin{abstract}
Pickands dependence functions characterize bivariate extreme value copulas.
In this paper, we study the class of polynomial Pickands functions. We
provide a
solution for the characterization of such polynomials of
degree at most $m+2$, $m\geq0$, and show that these can be
parameterized by a vector in $\mathbb{R}^{m+1}$ belonging to the intersection
of two ellipsoids. We also study the class of Bernstein approximations
of order $m+2$ of
Pickands functions which are shown to be (polynomial) Pickands
functions and
parameterized by a vector in
$\mathbb{R}^{m+1}$ belonging to a polytope. We give necessary and
sufficient conditions for which a polynomial Pickands
function is in fact a Bernstein approximation of some Pickands
function. Approximation results of Pickands dependence
functions by polynomials are given. Finally, inferential methodology is
discussed and comparisons based on simulated
data are provided.
\end{abstract}

%
\end{frontmatter}

\section{Introduction}
Bivariate extreme value copulas are characterized by the Pickands
dependence functions, these are
functions $A\dvtx [0,1] \to\reals$ which satisfy
the  conditions:
\begin{enumerate}[2.]
\item[1.] \textit{Boundary conditions}: $(1-t)\vee t\leq A(t) \leq1$,
$t\in[0,1]$.
\item[2.] \textit{Convexity condition}: $A$ is convex.
\end{enumerate}
In the presence of the \textit{convexity condition}, the above
\textit{boundary conditions} are simply saying that the lines $\ell_1(t)=1-t$
and $\ell_2(t)=t$, $t \in[0,1]$, are both \textit{support lines} of
$A$ at
$t=0$ and $t=1$ respectively, and it follows that they can be equivalently
replaced by

1. \textit{Endpoint conditions}: $A(0)=1=A(1)$ and $-1\leq A'(0)$,
$A'(1)\leq1$
(see, e.g., Roberts and Varberg \cite{RV73}, page 14, problem
K). Note also that in the case
where $A$ is twice differentiable the \textit{convexity condition} can be
replaced by $A''(t)\geq0$, for all $t\in[0,1]$.

Let $\mathcal{A}$ be the space of Pickands dependence functions. For
$A \in\mathcal{A}$, the copula is
%
%
\begin{equation}
\label{eq:CA} C_A(u,v)=\exp\biggl\{\log(uv)A \biggl(
\frac{\log v}{\log uv} \biggr) \biggr\}, \qquad0< u,v\leq1,
\end{equation}
and $C_A$ has a density if $A'$ is absolutely continuous on $[0,1]$.
Bivariate extreme value copulas may also be parameterized by the so-called
spectral measure $H$,
that is, a positive measure on $[0,1]$ such that
$\int_{[0,1]}w H(\mathrm{d}w)=1=\int_{[0,1]}(1-w) H(\mathrm{d}w)$.
Let $\mathcal{H}$ denote the space of spectral measures. The one-to-one
correspondence between $\mathcal{A}$ and $\mathcal{H}$, is given by
%
%
\begin{equation}
\label{bijection} A(t)=\int_{[0,1]}\bigl\{(1-t)w\bigr\} \vee\bigl
\{t(1-w)\bigr\} H(\mathrm{d}w),\qquad t\in[0,1],
\end{equation}
and $H([0,w])=1+A'_+(w)$, for all
$w\in[0,1)$, where $A'_+$ is the right derivative of $A$, see
Beirlant \textit{et~al.} \cite{BGST04}.

In the past literature, parametric models for the function $A$ have
been studied
for instance in Tawn \cite{Tawn88}, Coles and Tawn \cite{CT91},
Joe \textit{et al.} \cite{JSW92}, Ledford and Tawn \cite{LT96} and
Dupuis and Tawn \cite{DT01}. The modeling of $H$ has also been
considered, as it
offers some
advantages, especially when extensions to the
multivariate case are desired. Inference
within parametric families of spectral measures can be found in
Boldi and Davison \cite{BD07}, Coles and Tawn \cite{CT91},
Coles and Tawn \cite{CT94}, de~Haan \textit{et al.} \cite{dHNP08},
Einmahl \textit{et al.} \cite{EKS08},
Joe \textit{et al.} \cite{JSW92}, Ledford and Tawn \cite{LT96} and
Smith \cite{Smith94}. Polynomial splines
models have
also been proposed by Guillotte \textit{et al.} \cite{GPS11} and Foug\`{e}res \textit{et al.} \cite{FMN13}.

In view of the infinite-dimensional nature of the space of Pickands
functions, inference is often done
nonparametrically for maximum flexibility, see for instance Cap\'{e}ra\`{a} \textit{et al.} \cite{CFG97},
Hall and Tajvidi \cite{HT00}, Fils-Villetard \textit{et al.} \cite{FGS08} and
recently B{\"u}cher \textit{et al.} \cite{BDV11}, to cite only a few.
However, the nonparametric estimators do not usually satisfy the
properties of
Pickands functions for finite samples and are often modified to do so.
This can
be cumbersome, see, for instance, Fils-Villetard \textit{et al.} \cite{FGS08} and
the references therein.
One alternative to the nonparametric approaches is based on having a
series of nested
models indexed by $m$, $m\geq0$. Each model is parametric but $m$ is
not bounded. Within each model, a genuine
estimator of the Pickands dependence function is proposed. In this
spirit, we
consider modeling the Pickands function
$A$ using polynomials on $[0,1]$. The latter are natural candidates
here, the quadratic case
being known as the symmetric mixed model, while the cubic case
corresponds to
the asymmetric mixed model, see Beirlant \textit{et~al.} \cite{BGST04},
pages 308 and 309, and the references
therein. An attempt at characterizing the space of Pickands
polynomials for higher degrees was made in Kl{\"u}ppelberg and May \cite
{KM06}, page 1472, Theorem~2.5. It is stated there that a polynomial (in
the power basis) $A(t)=1-(\sum_{k=2}^ma_k)t+\sum_{k=2}^ma_kt^k$, $t
\in[0,1]$, is a Pickands dependence function if and
only if its coefficients satisfy the four following conditions:
\[
0\leq a_2,\qquad  0 \leq\sum_{k=2}^m
a_k,\qquad  0\leq\sum_{k=2}^m
(k-1)a_k\leq1 \quad \mbox{and}\quad  0\leq\sum_{k=2}^m
k(k-1)a_k.
\]
It turns out that these conditions fail to be sufficient for $m=4$. In
fact, a counterexample dates back to
Beirlant \textit{et~al.} \cite{BGST04}, page 308: the polynomial
$A(t)=1-t^3+t^4$ satisfies the above
conditions but is not convex. We will return to this example in
Section~\ref{sect:characterization}.

First, Section~\ref{sect:Bernstein} is a short preliminary section
giving the
essential definitions and general mathematical results concerning the Bernstein
basis used throughout the paper. The reason for this choice of basis
here is,
essentially, its ability to reflect geometric and
algebraic particularities of this problem on the coefficients.
This will be reinforced, with an example and pointers to the text in the
concluding remarks of Section~\ref{sect:conclusion}. In
Section~\ref{sect:characterization}, the aim is to characterize, for each
polynomial having degree at most $m$, $m\geq0$, the domain in which
the coefficients of the polynomial (in the
Bernstein basis) must
belong. We are looking for a constructive solution leading ultimately
to a
parameterization, which could be directly applicable in the search of the
maximum
likelihood estimate, for instance. We approach the problem by
exploiting the
one-to-one correspondence between the Pickands function and the related
spectral measure mentioned above. For polynomials, this boils down to
the link
between the
Pickands function and its second derivative (a nonnegative polynomial). The
solution therefore exploits the characterization of nonnegative polynomials
given by Luk\'acs, see Karlin and Shapley \cite{KS53}. We subsequently
refer to this
solution as
being the solution to the \textit{full model}. In
Section~\ref{sect:approximation}, the focus is on characterizing,
for each degree $m\geq0$, the space of polynomials that can be obtained
from a Bernstein approximation of a Pickands dependence function. We
call this
model the \textit{submodel}. We try to give an answer to the question:
what is the gap between the \textit{full model} and the \textit
{submodel}? Now
in Section~\ref{sect:quality}, we use tools from approximation theory and
probability for obtaining accurate bounds for measuring the closeness between
the space of Pickands dependence functions and the one obtained from the
\textit{submodel}. Note that Bernstein approximations of copulas have
appeared in the past literature and some of their properties have been studied
in Sancetta and Satchell \cite{SS01} and Sancetta and Satchell \cite{SS04}.
Finally, in Section~\ref
{sect:simulations}, we
present the results of a simulation comparing the maximum likelihood estimator
from the full model with the one from the submodel and a version of the
popular nonparametric
``CFG'' estimator in Cap{\'e}ra{\`a} \textit{et al.} \cite{CFG97}. Our estimators
should be easily implemented,
and could become part of standard extreme value statistics
packages, such as the EVD package in R. The existing packages offer only
polynomial models of degree at
most three, this may be explained by the absence of a successful
characterization for higher degrees in the previous
literature.

\section{Bernstein polynomials}
\label{sect:Bernstein}
Let $\mathcal{P}_m$ be the space of polynomials on the interval
$[0,1]$, with
degree at most $m$, $m\geq0$. We shall represent polynomials in
$\mathcal{P}_m$ using the commonly
called Bernstein basis. For self-containedness of the exposition, we
begin with
a short (preliminary) section describing some of the key properties of this
basis that will serve in the following developments. The connection between
this basis and the binomial distribution provides a very
powerful tool in some of the proofs. We are also going to exploit the
binomial identities given below.

\subsection{The basis}
Henceforth, let $S_n$ be a $\Binomial(n,x)$ random variable, $n\geq0$,
where the value of $x\in[0,1]$ will always
be clearly indicated by the context. The members of the Bernstein basis of
degree $m$, $m\geq0$, are the polynomials
\begin{eqnarray*}
b_{k,m}(x)={m\choose k}x^k(1-x)^{m-k}=
\prob_x(S_m=k),\qquad  x\in[0,1],
\end{eqnarray*}
for $0\leq k\leq m$. The coefficients of a polynomial $P\in\mathcal
{P}_m$ in
the Bernstein basis of degree $m$,
will be denoted by $\{c(k,m;P): k=0,\ldots,m\}$, and so for $P\in
\mathcal{P}_m$ we have
\[
P(x)=\sum_{k=0}^m c(k,m;P)b_{k,m}(x)=
\expec_x\bigl[c(S_m,m;P)\bigr],\qquad  x\in[0,1].
\]
Let $\Delta$ be the
forward difference operator. Here, when applied to a function $f$ of
two arguments: $(k,m) \mapsto f(k,m)$, it is
understood that $\Delta$ operates on the first argument:
$\Delta f(k,m)=f(k+1,m)-f(k,m)$ and
$\Delta^2 f(k,m)=f(k+2,m)-2f(k+1,m)+f(k,m)$, for all $(k,m)$.
Here are some basic results related to the Bernstein basis.
%

%
\begin{proposition}
\label{thm:Bern_prop}
For
$0\leq k\leq m$, we have
\begin{enumerate}[(iii)]
\item[(i)] $b'_{k,m}=m(b_{k-1,m-1}-b_{k,m-1})=-m\Delta
b_{k-1,m-1}$, for $m\geq1$, with the convention $b_{j,m}=0$ for $j
\notin
\{0,\ldots,m\}$,
\item[(ii)] for $0\leq\ell\leq n$, $b_{k,m}b_{\ell,n}=\frac{{m\choose
k}{n\choose
\ell}}{{m+n\choose k+\ell}}b_{k+\ell,m+n}$,
\item[(iii)]
$
\int_0^tb_{k,m}(w) \,\mathrm{d}w=\frac{1}{m+1}\sum
_{j=k+1}^{m+1}b_{j,m+1}(t)=\frac{1}{m+1}
-\int_t^1 b_{k,m}(w) \,\mathrm{d}w, t\in[0,1]$.
\end{enumerate}
\end{proposition}

The following are direct but useful consequences of the
above proposition. Properties (i) and (ii) (obtained using the summation
by parts formula) show how derivatives of polynomials represented in the
Bernstein basis act on the coefficients.
%

%
\begin{proposition}
\label{thm:Bern_calc}
Let $P\in\mathcal{P}_m$, with $P'$, $P''$, its first and second derivative
respectively.
\begin{enumerate}[(iii)]
\item[(i)] $c(k,m-1;P')=m\Delta c(k,m;P)$, $0\leq k\leq m-1$, $m\geq1$, so that
\[
P'(t)=m\expec_t\bigl\{\Delta c(S_{m-1},m;P)
\bigr\},\qquad  t\in[0,1],
\]
\item[(ii)] $c(k,m-2;P'')=m(m-1)\Delta^2 c(k,m;P)$, $0\leq k\leq m-2$,
$m\geq2$, so
that
\[
P''(t)=m(m-1)\expec_t\bigl\{
\Delta^2 c(S_{m-2},m;P)\bigr\},\qquad  t\in[0,1],
\]
\item[(iii)] for $0\leq i\leq m$ and $0\leq j\leq n$, and $t\in[0,1]$,
\begin{eqnarray*}
\int_0^t b_{i,m}(x)b_{j,n}(x)
\,\mathrm{d}x&=& \frac{{m\choose i}{n\choose j}}{{m+n\choose i+j}}
\frac{1}{m+n+1}\expec_t\bigl\{
\indic(i+j<S_{m+n+1})\bigr\},
\\
\int_t^1 b_{i,m}(x)b_{j,n}(x)
\,\mathrm{d}x &=& \frac{{m\choose i}{n\choose j}}{{m+n\choose i+j}}
\frac{1}{m+n+1}\expec_t\bigl\{
\indic(i+j\geq S_{m+n+1})\bigr\}.
\end{eqnarray*}
\end{enumerate}
\end{proposition}

\subsection{Bernstein approximations}
\label{sect:BernsteinApprox}
In Section~\ref{sect:approximation}, we will construct a submodel
based on
Bernstein approximations of Pickands functions. In general, the
$m$th-order Bernstein approximation of any function $f\dvtx [0,1]\to
\reals$
is given by
\[
B_m(f,t)=\sum_{k=0}^m f(k/m)
b_{k,m}(t)=\expec_t\bigl[f(S_m/m)\bigr], \qquad t
\in[0,1].
\]
It transforms $f$ into the polynomial
$B_m(f,\cdot)\in\mathcal{P}_m$, with coefficients equal to the
function values
on the uniformly spaced grid $\{k/m:k=0,\ldots,m\}$. Notice that $B_m(f,
0)=f(0)$ and $B_m(f, 1)=f(1)$. However, $B_m(f,\cdot)$ does not
necessarily interpolate $f$ on the grid. As a first implication of
Proposition~\ref{thm:Bern_calc} above, we get that the convexity of
$f$ implies
that of $B_m(f,\cdot)$.

\subsection{Binomial identities}
\label{sect:binident}
We will make extensive use of the following identities.
Let $X_{1},\ldots, X_{m}$ be independent $\Bernoulli(t)$ random
variables,
$t\in[0,1]$, and let $f\dvtx \{0,\ldots,m\}\to
\reals$ be a function. Here
$S_{m}=\sum_{k=1}^mX_{k}\sim\Binomial(m,t)$, and
%
%
\begin{eqnarray}
\label{eq:bin_expecfs} \expec_t\bigl\{f(S_{m})\bigr\}&=&
\expec_t\bigl\{f(S_{m-1}+X_{m})\bigr\} =
\expec_t\bigl[\expec_t\bigl\{f(S_{m-1}+X_{m})
\mid S_{m-1}\bigr\}\bigr]
\nonumber
\\[-8pt]\\[-8pt]
&=&t\expec_t\bigl\{f(S_{m-1}+1)\bigr\}+(1-t)
\expec_t\bigl\{f(S_{m-1})\bigr\},\nonumber
\\
\label{eq:bin_expecsfs} \expec_t\bigl\{S_{m}f(S_{m})
\bigr\}&=&m\expec_t\bigl\{X_{m}f(S_{m-1}+X_{m})
\bigr\} =m\expec_t\bigl[\expec_t\bigl
\{X_{m}f(S_{m-1}+X_{m})\mid S_{m-1}\bigr\}
\bigr]
\nonumber
\\[-8pt]\\[-8pt]
&=&mt\expec_t\bigl\{f(S_{m-1}+1)\bigr\},\nonumber
\end{eqnarray}
so that from \eqref{eq:bin_expecfs} and \eqref{eq:bin_expecsfs}, we have
$\expec_t\{(m-S_{m})f(S_{m})\}=m(1-t)\expec_t\{f(S_{m-1})\}$.
These identities are special cases of more general
identities developed for the exponential families, see Hudson \cite{H78}.

\section{The characterization}
\label{sect:characterization}
The problem of characterizing the polynomial Pickands functions is
solved in
Section~\ref{sect:PolPick} below. To do so, we will take further
advantage of
the one-to-one correspondence between the spectral measure $H$ and the
related Pickands
function $A$. In fact, under absolute continuity of $A'$, we will show
that the
problem of modeling $A$ (or $H$) boils down to modeling some nonnegative
function $h$
on $(0,1)$ satisfying the condition
%
%
\begin{equation}
\label{function_h} \biggl\{\int_0^1 (1-w) h(w)\,
\mathrm{d}w \biggr\}\vee\biggl\{\int_0^1 w
h(w)\,\mathrm{d}w \biggr\}\leq1.
\end{equation}
Although this is (briefly) mentioned in Beirlant \textit{et~al.} \cite
{BGST04}, page 269, we
provide the
details in Section~\ref{sect:representation} because the use of this special
function $h$ is a key element in the paper. Note that this framework is
slightly more general then what is really needed for obtaining the
characterization later on, because we will be working with polynomials rather
than merely absolutely continuous functions. However, it is effortless
to do so
here.

\subsection{A representation theorem}
\label{sect:representation}
Let $A \in\mathcal{ A}$, let $H$ on $[0,1]$ be the spectral measure
associated with $A$ and let $\mu$ be the Lebesgue measure on $[0,1]$. Assume
that $A'$ is absolutely continuous on $[0,1]$, then, in particular, the copula
$C_A$ in equation (\ref{eq:CA}) has a density with respect to the Lebesgue
measure on $[0,1]^2$. In this case, we derive a convenient integral
representation for $A(t)$, $t\in[0,1]$, in terms of the Radon--Nikod\'ym
derivative $h$ of the restriction of $H$ on $(0,1)$. Note that $H$ may still
have point masses at 0 and 1, which can be seen by $H(\{0\})=1+A'(0)$ and
$H(\{1\})=1-A'(1)$. Under
the above regularity condition on $A'$, we will show that the knowledge (almost
everywhere) of $A''$ alone is enough to evaluate $A$, and similarly,
the knowledge of the Radon--Nikod\'ym derivative $h$ of the restriction
of $H$ on
$(0,1)$ is enough to know the measure $H$. In particular, we have $h=A''$
almost everywhere on $(0,1)$.

%
%
\begin{theorem}[(Integral representation)]
\label{thm:representation}
\begin{enumerate}[(ii)]
\item[(i)] Let $A\in\mathcal{A}$, let $H$ be the spectral measure
associated with $A$. Assume that $A'$ is absolutely continuous and let
$h=A''$ almost everywhere on $(0,1)$. We have
%
%
\begin{equation}
\label{A_representation} A(t)= 1-\int_0^1\bigl[\bigl
\{(1-t)w\bigr\}\wedge\bigl\{t(1-w)\bigr\}\bigr] h(w) \,\mathrm
{d}w,\qquad  t\in[0,1],
\end{equation}
and the function $h$ satisfies the
condition \eqref{function_h}. Let $\mu$ be the Lebesgue measure and
$\delta_x$
denote a Dirac measure at $x$, the spectral measure is given by
$H=h_0\delta_0+\mathring{H}+h_1\delta_1$, with
%
%
\begin{eqnarray}
\label{H_representation} %
\mathrm{d}\mathring{H}/\mathrm{d}\mu&=& h,\nonumber
\\
h_0&=&H\bigl(\{0\}\bigr)=1-\int_0^1
(1-w) h(w) \,\mathrm{d}w,
\\
h_1&=&H\bigl(\{1\}\bigr)=1-\int_0^1
w h(w) \,\mathrm{d}w.\nonumber %
\end{eqnarray}
\item[(ii)] Conversely, let $h\dvtx (0,1)\to\reals$
be a nonnegative function satisfying condition \eqref{function_h}. The measure
$H$ given
by \eqref{H_representation} is a spectral measure, its corresponding Pickands
function $A$ is given by \eqref{A_representation}, and we have
$A''=h$ almost everywhere on $(0,1)$.
\end{enumerate}
\end{theorem}

\begin{pf}
(i) Integration by parts yields
%
%
\begin{equation}
\label{0tot} \int_0^t w h(w) \,\mathrm{d}w
=tA'(t)+\bigl\{1-A(t)\bigr\},\qquad  t\in[0,1],
\end{equation}
and
%
%
\begin{equation}
\label{tto1} \int_t^1 (1-w) h(w) \,
\mathrm{d}w =\bigl\{1-A(t)\bigr\}-(1-t)A'(t),\qquad  t\in[0,1],
\end{equation}
and by combining these two equations, we obtain
\begin{eqnarray*}
\bigl\{1-A(t)\bigr\} &=& (1-t)\int_0^t w h(w)
\,\mathrm{d}w +t\int_t^1 (1-w) h(w) \,
\mathrm{d}w
\\
&=& \int_0^1\bigl\{(1-t)w\bigr\}\wedge\bigl
\{t(1-w)\bigr\}h(w)\,\mathrm{d}w,\qquad  t\in[0,1].
\end{eqnarray*}
Letting $t=1$ in \eqref{0tot} and $t=0$ in \eqref{tto1} shows that the
condition \eqref{function_h} is satisfied:
\[
\int_0^1w h(w)=A'(1)\leq1
\quad \mbox{and}\quad  \int_0^1(1-w) h(w)=-A'(0)
\leq1.
\]
From the one-to-one correspondence between Pickands
functions and spectral measures, the only thing left to show is that if
$H$ is
the measure given by \eqref{H_representation}, then it is the spectral measure
associated to $A$, that is, equation \eqref{bijection} holds:
%
%
\begin{eqnarray}
\label{identity} A(t) &=& 1-\int_0^1\bigl[\bigl
\{(1-t)w\bigr\}\wedge\bigl\{t(1-w)\bigr\}\bigr] h(w) \,\mathrm{d}w
\nonumber
\\
&=& (1-t) \biggl(1-\int_0^tw h(w) \,
\mathrm{d}w \biggr)+t \biggl(1-\int_t^1(1-w)
h(w) \,\mathrm{d}w \biggr)
\nonumber
\\
&=& (1-t) \biggl(h_1+\int_t^1w
h(w) \,\mathrm{d}w \biggr)+t \biggl(h_0+\int_0^t(1-w)
h(w) \,\mathrm{d}w \biggr)
\\
&=& (1-t)h_1+th_0+\int_0^1
\bigl[\bigl\{(1-t)w\bigr\}\vee\bigl\{t(1-w)\bigr\}\bigr] h(w) \,
\mathrm{d}w
\nonumber
\\
&=& \int_{[0,1]}\bigl\{(1-t)w\bigr\} \vee\bigl\{t(1-w)\bigr\}
H(\mathrm{d}w), \qquad t\in[0,1],\nonumber
\end{eqnarray}
as claimed.

(ii) For the converse, if $h\dvtx (0,1)\to[0,\infty)$
satisfies the condition \eqref{function_h}, then the measure $H$ given by
\eqref{H_representation} is a spectral
measure: it is clearly positive, and
\[
\int_{[0,1]} w H(\mathrm{d}w)=h_1+\int
_0^1wh(w)\,\mathrm{d}w=1 =h_0+\int
_0^1(1-w)h(w)\,\mathrm{d}w =\int
_{[0,1]} (1-w) H(\mathrm{d}w). %
\]
Let $A$ be its corresponding Pickands function given by
equation \eqref{bijection}. The above equalities that lead to \eqref{identity}
show that $A$ is also given by \eqref{A_representation}. Finally, a direct
calculation shows that $A''=h$ almost everywhere on $(0,1)$.
\end{pf}

\subsection{The space of polynomial Pickands functions}
\label{sect:PolPick}
The space $\mathcal{A}_m=\mathcal{P}_m\cap\mathcal{A}$
corresponds to polynomial Pickands functions with degree at most $m$,
$m\geq
0$. The reader should refer to Section~\ref{sect:Bernstein} for
notations used
in the following. It is
immediate that $\mathcal{A}_0=\mathcal{A}_1$, and $A\in\mathcal{A}_0$
if and only if $A(t)=1$, for all $t\in[0,1]$.
The nontrivial cases start with Pickands polynomials of degree at least two.
We will characterize $\mathcal{A}_m$, for all $m\geq2$.
Note that in the Bernstein basis, the \textit{endpoint conditions} are
quite simple to verify because they are directly related to
the first and last two coefficients only.
%

%
\begin{proposition}
\label{thm:pol_constraints}
Let $A\in\mathcal{P}_{m+2}$, $m\geq0$. We have $A\in\mathcal{A}_{m+2}$
if and only if the two following conditions are verified:
\begin{enumerate}[2.]
\item[1.] \textit{Endpoint conditions}:
\[
c(0,m+2;A)=1=c(m+2, m+2;A),
\]
and
\[
c(1,m+2;A)\wedge c(m+1,m+2;A)\geq(m+1)/(m+2).
\]
\item[2.] \textit{Convexity condition}: $A''(t)\geq0$, $t\in[0,1]$.
\end{enumerate}
\end{proposition}

\begin{pf}
We need to verify the \textit{endpoint conditions}, that is, $A(0),A(1)=1$,
$-A'(0),A'(1)\leq1$.
In the Bernstein basis, by Proposition~\ref{thm:Bern_calc}, we have
\begin{eqnarray*}
c(0,m+2;A)&=&A(0)=1 \quad \mbox{and}\quad  c(m+2, m+2;A)=A(1)=1,
\\
1&\geq&-A'(0)=-c\bigl(0,m+1;A'\bigr)=-(m+2)\Delta
c(0,m+2;A)\\
&=&(m+2)\bigl\{1-c(1,m+2;A)\bigr\},
\end{eqnarray*}
and
\begin{eqnarray*}
1&\geq& A'(1)=c\bigl(m+1,m+1;A'\bigr)=(m+2)\Delta
c(m+1,m+2;A)\\
&=&(m+2)\bigl\{ 1-c(m+1,m+2;A)\bigr\}.
\end{eqnarray*}
\upqed
\end{pf}

It is therefore the \textit{convexity condition} above that is more delicate
to express in terms of the coefficients of $A$. However, the
spaces $\mathcal{A}_2$ and $\mathcal{A}_3$ are still easy to
characterize, as
can be seen in the following example.
%
%

%
\begin{Example}[($\boldsymbol{\mathcal{A}_2}$ and $\boldsymbol{\mathcal{A}_3}$)]
\label{ex:A2A3}
We have $A\in
\mathcal{A}_2$ if and only if
$c(0,2;A)=1=c(2,2;A)$ and $1/2\leq c(1,2;A)\leq1$. To see this,
let $A\in\mathcal{P}_2$. The polynomial $A$ satisfies the \textit{endpoint
conditions}
of Proposition~\ref{thm:pol_constraints} if and only if
$c(0,2;A)=1=c(2,2;A)$ and $1/2\leq c(1,2;A)$. Moreover $A''(t)=A''(0)$, for
all $t\in[0,1]$, and by Proposition~\ref{thm:Bern_calc}, $
A''(0)=2\Delta^2
c(0,2;A)=4\{1-c(1,2;A)\}$. Therefore, $A$ is convex if and only if
$c(1,2;A)\leq1$.

We have $A\in\mathcal{A}_3$ if and only if
$c(0,3;A)=1=c(3,3;A)$, and the couple $(c(1,3;A),\allowbreak c(2,3;A))$ belongs to the
polytope derived from the four linear inequalities:
\[
c(1,3;A) \wedge c(2,3;A)\geq2/3
\]
and
\[
\bigl\{2c(1,3;A)-c(2,3;A)\bigr\} \vee\bigl\{2c(2,3;A)-c(1,3;A)\bigr
\}\leq1.
\]
To show this, let $A\in\mathcal{P}_3$. The polynomial $A$ satisfies the
\textit{endpoint conditions} of Proposition~\ref{thm:pol_constraints}
if and
only if
$c(0,3;A)=1=c(3,3;A)$, and $c(1,3;A) \wedge c(2,3;A)\geq2/3$. Also, since
$A''(t)=(1-t)A''(0)+tA''(1)$, $t\in[0,1]$, we have $A''\geq0$, in $[0,1]$
if and only if $A''(0)\wedge A''(1)\geq0$. By evaluating
$A''(k)=6\Delta^2
c(k,3;A)$ for $k=0,1$, we get that $A$ is convex
if and only if $\{2c(1,3;A)-c(2,3;A)\} \vee\{2c(2,3;A)-c(1,3;A)\}\leq1$.
\end{Example}

For $m\in\{0,1\}$, finding $\mathcal{A}_{m+2}$ explicitly (in
Example~\ref{ex:A2A3} above) is easy mostly because
$A''\in\mathcal{P}_1$ is nonnegative if and only if
$A''(0) \wedge A''(1)\geq0$. When $A\in\mathcal{A}_{m+2}$ with $m>1$ things
become more complicated, as the following example illustrates.
%
%

%
\begin{Example}[(The counterexample polynomial in $\boldsymbol{\mathcal{P}_4}$)]
\label{ex:counterexPol}
Consider the polynomial
$A(t)=1-t^3+t^4$, which served as the counterexample provided
by Beirlant \textit{et~al.} \cite{BGST04} discussed in the
introduction. Here,
$A''(t)=12t(t-1/2)$, so $A''(0)\wedge A''(1)=0$ but
$\{t\colon A''(t)\geq0\}=\{0\}\cup[1/2,1]\neq[0,1]$, and therefore
$A\in\mathcal{P}_4\setminus\mathcal{A}_4$.
\end{Example}

Essentially, it all boils down to finding the coefficients making $A''$
nonnegative, and the problem is now reoriented towards the
characterization of
such $A''$. We will essentially be using the \textit{Integral representation}
Theorem~\ref{thm:representation} and a characterization of
the nonnegative polynomials on $[0,1]$ due to Luk\'{a}cs. The main result
is given in Theorem~\ref{thm:mainthm}, while a
parametric representation of the polynomial Pickands functions and the geometric
shape of the corresponding parameter space is provided via
Theorem~\ref{thm:parameterization}. To get to these results, we first need
two technical lemmas, namely Lemma~\ref{thm:HtoA} and Lemma~\ref{thm:PandQ}.

Because of the new orientation of the problem as mentioned above, it is
convenient to introduce new spaces $\mathcal{H}_m\subset\mathcal{P}_m$,
$m\geq0$. We say that a polynomial $h$ belongs to
$\mathcal{H}_m$, $m\geq0$, when the  following conditions are satisfied:
\begin{enumerate}[2.]
\item[1.]\textit{Endpoint derivatives conditions:}
%
%
\begin{eqnarray}
\label{eq:boundcondh} \frac{1}{m+1} \sum_{j=0}^m
\biggl(1-\frac{j+1}{m+2} \biggr)c(j,m;h)&\leq&1 \quad \mbox{and}\quad
 \frac{1}{m+1}
\sum_{j=0}^m \frac{j+1}{m+2} c(j,m;h)
\leq1.\qquad
\end{eqnarray}
\item[2.]\textit{Nonnegativity condition:} $h(t)\geq0$, $t\in[0,1]$.
\end{enumerate}

\begin{Remark*} These are the same conditions as those given at the very
beginning of this section when restricted to polynomials, in
particular, the
\textit{endpoint derivatives conditions} \eqref{eq:boundcondh}
come directly from the condition \eqref{function_h}, by using
Proposition~\ref{thm:Bern_prop}. When a polynomial $A$ on $[0,1]$
has an integral reprentation \eqref{A_representation}, the
condition \eqref{function_h} corresponds to
$-1\leq A'(0), A'(1)\leq1$, hence the name \textit{endpoint derivatives
conditions}.
From the expression \eqref{A_representation}
we get the following result (notice the resemblance between the
polynomial $A$
and its coefficients in the Bernstein basis\ldots).
\end{Remark*}

%
\begin{lemma}
\label{thm:HtoA}
Let $m\geq0$.
\begin{enumerate}[(ii)]
\item[(i)] If $h\in\mathcal{P}_{m}$ and
%
%
\begin{equation}
\label{eq:PolA} A(t)= 1-\int_0^1\bigl[\bigl
\{(1-t)w\bigr\}\wedge\bigl\{t(1-w)\bigr\}\bigr] h(w) \,\mathrm
{d}w,\qquad  t\in[0,1],
\end{equation}
then $A\in\mathcal{P}_{m+2}$ and its $k$th coefficient,
$k=0,1,\ldots,m+2$,
in the Bernstein basis is
%
%
\begin{eqnarray}
\label{eq:HtoA} &&c(k,m+2;A)\nonumber\\[-8pt]\\[-8pt]
&&\quad  =1-\frac{1}{m+1}\sum_{j=0}^m
\biggl[ \biggl\{ \biggl(1-\frac{k}{m+2} \biggr)\frac{j+1}{m+2}
\biggr\}
\wedge\biggl\{\frac{k}{m+2} \biggl(1-\frac{j+1}{m+2} \biggr)
\biggr\}
\biggr]c(j,m;h).\nonumber
\end{eqnarray}

\item[(ii)] If $A\in\mathcal{P}_{m+2}$ and $h=A''$, then $h\in\mathcal
{P}_{m}$ and
\[
c(k,m;h)=(m+2) (m+1)\Delta^2 c(k,m+2;A), \qquad k=0,\ldots,m.
\]
\end{enumerate}
\end{lemma}

\begin{pf}
(i)
We have the equality
%
%
\begin{eqnarray}
\label{eq:inteq} &&\int_0^1\bigl[\bigl\{(1-t)w
\bigr\}\wedge\bigl\{t(1-w)\bigr\}\bigr] h(w) \,\mathrm{d}w\nonumber\\[-8pt]\\[-8pt]
&&\quad = \int
_0^tw h(w) \,\mathrm{d}w+ t\int
_t^1 h(w) \,\mathrm{d}w- t\int
_0^1w h(w) \,\mathrm{d}w.\nonumber
\end{eqnarray}
Using Proposition~\ref{thm:Bern_calc}(iii) and the binomial
identities \eqref{eq:bin_expecsfs},
we obtain:
\begin{eqnarray*}
\int_0^tw h(w) \,\mathrm{d}w &=&
\frac{1}{m+1}\sum_{j=0}^m
\frac{j+1}{m+2}\expec_t \bigl\{\indic(j+1< S_{m+2}) \bigr
\}c(j,m;h),
\\
t\int_t^1 h(w) \,\mathrm{d}w &=&
\frac{t}{m+1}\sum_{j=0}^m
\expec_t\bigl\{\indic(j\geq S_{m+1})\bigr\} c(j,m;h)
\\
&=& \frac{1}{m+1}\sum_{j=0}^m
\expec_t \biggl\{\frac{S_{m+2}}{m+2} \indic(j+1\geq S_{m+2})
\biggr\}c(j,m;h),
\\
t\int_0^1w h(w) \,\mathrm{d}w &=&
\frac{1}{m+1}\sum_{j=0}^m
\frac{j+1}{m+2}\expec_t \biggl(\frac{S_{m+2}}{m+2} \biggr)c(j,m;h),
\end{eqnarray*}
and putting the three pieces together shows that \eqref{eq:inteq} equals
\[
\expec_t \Biggl(\frac{1}{m+1}\sum_{j=0}^m
\biggl[ \biggl\{ \biggl(1-\frac{S_{m+2}}{m+2} \biggr)\frac
{j+1}{m+2} \biggr\}
\wedge\biggl\{\frac{S_{m+2}} {m+2} \biggl(1-\frac{j+1}{m+2} \biggr
) \biggr\}
\biggr]c(j,m;h) \Biggr),
\]
which in turn equals $\expec_t\{1-c(S_{m+2},m+2;A)\}$, and
this gives \eqref{eq:HtoA}.

(ii) This is directly obtained from Proposition~\ref{thm:Bern_calc}.
\end{pf}

%
%
\begin{Example}[($\boldsymbol{\mathcal{A}_4}$)]
\label{ex:A4}
We can show directly that $h\in\mathcal{P}_2$ is
nonnegative on $[0,1]$ if and only if
\[
c(0,2;h)\wedge c(2,2;h)\geq0 \quad \mbox{and} \quad c(1,2;h)\geq-\sqrt
{c(0,2;h)c(2,2;h)}.
\]
Indeed, if
$h\in\mathcal{P}_2$ is
nonnegative on $[0,1]$, then necessarily $c(0,2;h)=h(0)\geq0$ and
$c(2,2;h)=h(1)\geq0$. Since the latter conditions need to be satisfied,
lets assume them and look at the remaining coefficient $c(1,2;h)$. If
$c(1,2;h)\geq0$, then $h$ is nonnegative on $[0,1]$. If however $c(1,2;h)<0$,
then we have $\Delta^2c(0,2;h)>0$ and $-\Delta
c(0,2;h)/\Delta^2c(0,2;h)\in(0,1)$. By
expressing
\[
h(t)=\frac{c(0,2;h)c(2,2;h)-c(1,2;h)^2}{\Delta^2c(0,2;h)}+\Delta
^2c(0, 2;h) \biggl(t+
\frac{\Delta c(0,2;h)}{\Delta^2c(0,2;h)} \biggr)^2,\qquad  t\in[0,1],
\]
we see that $h$ is convex and attains its minimum on $(0,1)$. The result
follows from the sign of the minimum. The above conditions on
$c(0,2;h), c(1,2;h)$ and $c(2,2;h)$, together with the \textit{endpoint
derivatives conditions} and a direct application of
Lemma~\ref{thm:HtoA} gives the characterization of $\mathcal{A}_4$.
\end{Example}

Notice that, so far (for $m\leq4$), we have
characterized $\mathcal{A}_m$ using only elementary
mathematics. For higher degrees however, we use a
finer result by  Luk\'{a}cs,
concerning the characterization of nonnegative polynomials, which can
be found
in for instance Karlin and Shapley \cite{KS53}, Szeg{\H{o}}
\cite{Szego75} and P{\'o}lya and Szeg{\H{o}} \cite{PS98}. It says
that a
polynomial $h\in\mathcal{P}_m$ of degree at most $m$,
is nonnegative if and only if there are polynomials $P$ and $Q$
with $\deg(P)\leq\lfloor m/2 \rfloor$ and $\deg(Q)\leq\lfloor(m-1)/2
\rfloor$, such that
%
%
\begin{equation}
\label{eq:h} h(t)= %
\cases{ P^2(t)+t(1-t)Q^2(t),
& \quad \mbox{if $m$ is even},
\cr
tP^2(t)+(1-t)Q^2(t), &\quad \mbox{if $m$ is odd}, } %
\end{equation}
for all $t\in[0,1]$. The following lemma exploits Luk\'{a}cs' result and
provides expressions for the coefficients, in the Bernstein basis, of
nonnegative $h\in\mathcal{P}_{m}$,
$m\geq1$, using the Hypergeometric distribution. Here $Y\sim
\Hypergeo(n,M,N)$ if
\[
P(Y=k)=\frac{{M\choose k}{N-M\choose n-k}}{{N\choose n}}, \qquad (n+M-N)\vee
0\leq k \leq n\wedge M, 0\leq n,M\leq N.
\]
%

%
\begin{lemma}
\label{thm:PandQ}
A polynomial $h\in\mathcal{P}_{m}$, $m\geq1$, is nonnegative on
$[0,1]$ if
and only if there exist polynomials $P$ and $Q$ such that, the
coefficients of
$h$ are related to the coefficients of $P$ and $Q$ (in the Bernstein
basis) in
the following way: (for notational simplicity here, let $h(k,m)=c(k,m;h)$,
$p(k,m)=c(k,m;P)$ and $q(k,m)=c(k,m;Q)$ denote the coefficients)
\begin{enumerate}[(ii)]
\item[(i)] when $m$ is even, $P\in\mathcal{P}_{m/2}$, $Q\in\mathcal{P}_{(m-2)/2}$
and
%
%
\begin{eqnarray}
\label{eq:even} h(k,m)&=&\expec\biggl\{p \biggl(Y_1,\frac{m}{2}
\biggr)p \biggl(k-Y_1,\frac{m}{2} \biggr)\nonumber\\[-8pt]\\[-8pt]
&&\hphantom{\expec\biggl\{}{}+ \frac{k(m-k)}{m(m-1)}q
\biggl(Y_2,\frac{m-2}{2} \biggr)q \biggl(k-Y_2-1,
\frac{m-2}{2} \biggr) \biggr\},\nonumber
\end{eqnarray}
with
$Y_1\sim\Hypergeo(k,m/2,m)$, $k=0,\ldots,m$, and
$Y_2\sim\Hypergeo(k-1,(m-2)/2,m-2)$,
$k=1,\ldots,m-1$,
\item[(ii)] when $m$ is odd, $P,Q\in\mathcal{P}_{(m-1)/2}$
and
%
%
\begin{eqnarray}
\label{eq:odd} h(k,m)&=&\frac{1}{m}\expec\biggl\{kp \biggl(Y_1,
\frac{m-1}{2} \biggr)p \biggl(k-Y_1-1, \frac{
m-1}{2} \biggr)\nonumber\\[-8pt]\\[-8pt]
&&\hphantom{\frac{1}{m}\expec\biggl\{}{}+ (m-k)q \biggl(Y_2,\frac{m-1}{2} \biggr)q \biggl(k-Y_2,
\frac
{m-1}{2} \biggr) \biggr\},\nonumber
\end{eqnarray}
with
$Y_1\sim\Hypergeo(k-1,(m-1)/2,m-1)$, $k=1,\ldots,m$, and
$Y_2\sim\Hypergeo(k,(m-1)/2,m-1)$,
$k=0,\ldots,m-1$.
\end{enumerate}
\end{lemma}

\begin{pf}
The proof is a bit technical, and can be skipped without
affecting the readability of what follows. See the \hyperref[app]{Appendix}.
\end{pf}

By putting everything together, we get our main theorem.
%

%
\begin{theorem}[(Characterization)]
\label{thm:mainthm}
\begin{enumerate}[(iii)]
\item[(i)] For $m\geq1$, we have $h\in\mathcal{H}_m$ if and only if $h$ satisfies
the endpoint derivatives conditions \eqref{eq:boundcondh}
and there exist polynomials $P$ and $Q$ such that the coefficients of
$h$ (in
the Bernstein basis)
are given by \eqref{eq:even} when $m$ is even or by \eqref{eq:odd}
when $m$ is odd.
\item[(ii)] If $h\in\mathcal{H}_m$, $m\geq0$, then $A$ given by \eqref
{eq:PolA} is
in $\mathcal{A}_{m+2}$. Lemma~\ref{thm:HtoA} gives the coefficients
of $A$ in
terms of the  ones of $h$.
\item[(iii)] If $A\in\mathcal{A}_{m+2}$, $m\geq0$, then $A''\in\mathcal{H}_m$.
Lemma~\ref{thm:HtoA} gives the coefficients of $h$ in terms of the
ones of
$A$.
\end{enumerate}
\end{theorem}

A parameterization of $\mathcal{H}_m$ (and therefore of the polynomial Pickands
functions) can therefore be made via the coefficients of the
polynomials $P$ and $Q$
above. The following theorem describes the corresponding parameter space.
%

%
\begin{theorem}[(Parameterization)]
\label{thm:parameterization}
Let $m\geq0$, for polynomials $P\in\mathcal{P}_{\lfloor m/2 \rfloor
}$ and when $m\geq1$
$Q\in\mathcal{P}_{\lfloor(m-1)/2 \rfloor}$,
let $\theta$ be the concatenation of the coefficients of $P$ and the
ones of
$Q$ in the Bernstein basis, that is,
\[
\theta(k)= %
\cases{ c\bigl(k-1,\lfloor m/2 \rfloor;P\bigr), & \quad \mbox{if $ 1\leq k \leq\lfloor m/2 \rfloor+1$},
\cr \vspace{-8pt}
\cr
c\bigl(k-\lfloor m/2 \rfloor-2,\bigl
\lfloor(m-1)/2 \bigr\rfloor;Q\bigr), &\quad \mbox{if $ \lfloor m/2
\rfloor+2\leq k \leq
m+1$, $m\geq1$}. } %
\]
Let $h_{\theta}\in\mathcal{P}_{m}$ be constructed using polynomials
$P$ and $Q$
in formula \eqref{eq:h}. When $m\geq
1$, the coefficients of $h_{\theta}$ are given in Lemma~\ref
{thm:PandQ}. The parameter space
$\Theta_m=\{\theta\colon
h_{\theta}\in\mathcal{H}_m\}$ is given by
%
%
\begin{equation}
\label{eq:Thetam} \Theta_m= E_0 \cap E_1,
\end{equation}
where $E_0$ and $E_1$ are two ellipsoids in $\mathbb{R}^{m+1}$.
\end{theorem}

\begin{pf}
For every fixed value of $t\in[0,1]$, the function
$\theta\mapsto h_{\theta}(t)$ is a positive semidefinite quadratic form.
In fact, using the Bernstein basis and \eqref{eq:h},
we find that $h_{\theta}(t)$ is given by
\[
\cases{ \Biggl(\displaystyle \sum_{k=1}^{(m/2)+1}
b_{k-1,m/2}(t)\theta_k \Biggr)^2\cr
\quad {}+ \Biggl(\displaystyle \sum
_{k=1}^{m/2} \sqrt{t(1-t)}b_{k-1,(m/2)-1}(t)
\theta_{k+(m/2)+1} \Biggr)^2, & \quad \mbox{if $m$ is even,}
\cr
\Biggl(
\displaystyle \sum_{k=1}^{(m+1)/2} \sqrt{t}b_{k-1,(m-1)/2}(t)
\theta_k \Biggr)^2\cr
\quad {}+ \Biggl(\displaystyle \sum
_{k=1}^{(m+1)/2} \sqrt{1-t}b_{k-1,(m-1)/2}(t)
\theta_{k+(m+1)/2} \Biggr)^2, &\quad \mbox{if $m$ is odd.} } %
\]
Note that $h_{\theta}$ is the zero polynomial if and only if $\theta
=0$, and so
if $\|\theta\|>0$, then the function $t\mapsto
h_{\theta}(t)$ is positive except perhaps at a finite number of points
(roots). Now let $A_{\theta}$ be the polynomial
obtained from formula \eqref{eq:HtoA} with $h$ being replaced by
$h_{\theta}$. Since
\[
-A'_{\theta}(0)=\int_0^1(1-w)h_{\theta}(w)
\,\mathrm{d}w \quad \mbox{and}\quad  A'_{\theta}(1)=\int
_0^1wh_{\theta}(w) \,\mathrm{d}w,
\]
if $\|\theta\|>0$, it follows that $-A'_{\theta}(0)$ and $A'_{\theta
}(1)$ are both positive definite quadratic forms.
Finally, Theorem~\ref{thm:mainthm} says that $h_{\theta}\in\mathcal
{H}_m$ if and
only if the \textit{endpoint derivatives conditions} are verified:
\[
\biggl\{\int_0^1(1-w)h_{\theta}(w) \,
\mathrm{d}w \biggr\} \vee\biggl\{\int_0^1wh_{\theta}(w)
\,\mathrm{d}w \biggr\}\leq1,
\]
and so the sets
\[
E_0=\bigl\{\theta\in\reals^{m+1}\colon-A'_{\theta}(0)
\leq1\bigr\} \quad \mbox{and}\quad  E_1=\bigl\{\theta\in\reals^{m+1}
\colon A'_{\theta}(1)\leq1\bigr\},
\]
are therefore both ellipsoids in $\mathbb{R}^{m+1}$ and $\Theta_m=
E_0 \cap
E_1$.
\end{pf}

\section{The Bernstein approximations submodel}
\label{sect:approximation}
Let
\[
\mathcal{A}_{m}^+=\bigl\{B_{m}(A,\cdot): A\in\mathcal{A}
\bigr\},\qquad  m\geq0,
\]
be the set of Bernstein's approximations of Pickands functions. The
second part
of the following lemma gives a useful necessary and sufficient geometric
condition on the coefficients (in the Bernstein basis) for a polynomial
$A\in
\mathcal{P}_m$ to
belong to $\mathcal{A}_{m}^+$. The first part is used for this characterization
and is also useful in the next propositions.
%

%
\begin{lemma}
\label{thm:PiecewisePick}
For any polynomial $A\in\mathcal{P}_m$, $m \geq1$, represented in
the Bernstein
basis by $A(t)=\expec_t\{c(S_m,m;A)\}$, $t\in[0,1]$, let $A^*$ be the
piecewise linear function interpolating the points $\{(k/m,c(k,m;A))\}_{k=0}^m$,
\[
A^*(t)=\bigl(\lfloor mt \rfloor+1-mt\bigr)c\bigl(\lfloor mt \rfloor
,m;A\bigr)+
\bigl(mt-\lfloor mt \rfloor\bigr)c\bigl(\lfloor mt \rfloor
+1,m;A\bigr),\qquad  t\in[0,1].
\]
\vspace*{-\baselineskip}%
\begin{enumerate}[(ii)]
\item[(i)] If $\Delta^2c(k,m;A)\geq0$, $k=0,\ldots,m-2$, $m\geq2$, then
$A^*$ is
convex,
\item[(ii)](characterization) $A\in\mathcal{A}_{m}^+$ if and only if
$A^*\in\mathcal{A}$.
\end{enumerate}
\end{lemma}

\begin{pf}
(i)
 Consider the step function $\varphi(t)=m\Delta c(\lfloor mt
\rfloor\wedge(m-1),m;A)$, $t\in[0,1]$. Geometrically, $\varphi$ corresponds
to the right derivative of $A^*$ on [0,1). Now, $\Delta^2c(k,m;A)\geq0$,
$k=0,\ldots,m-2$, implies that $k\mapsto\Delta c(k,m;A)$,
$k=0,\ldots,m-1$, is nondecreasing, and therefore $\varphi$ is
nondecreasing on
$[0,1]$. Since $A^*(t)=1+\int_0^t \varphi(x) \,\mathrm{d}x$, for all
$t\in
[0,1]$, it
follows that $A^*$ is convex.

(ii) We have $A(\cdot)=B_m(A^*,\cdot)$. If $A^*\in\mathcal{A}$,
then $B_m(A^*,\cdot)\in\mathcal{A}_{m}^+$ and therefore $A\in
\mathcal{A}_{m}^+$. For
the converse, if $A \in\mathcal{A}_{m}^+$, then there is some $A_0
\in\mathcal{A}$ such that $A(\cdot)=B_m(A_0,\cdot)$. For $m=1$, the
fact that
$A^*\in\mathcal{A}$ follows trivially. For $m\geq2$, the convexity
of $A_0$
implies
\[
\Delta^2 c(k,m;A)=A_0 \biggl(\frac{k+2}{m}
\biggr)-2A_0 \biggl(\frac{k+1}{m} \biggr) +A_0
\biggl(\frac{k}{m} \biggr) \geq0,\qquad  k=0,\ldots,m-2,
\]
and (i) above shows that $A^*$ is convex. Also,
\[
A^*(0)=c(0,m;A)=A_0(0)=1=A_0(1)=c(m,m;A)=A^*(1)=1,
\]
and finally
\[
\bigl(A^*\bigr)'(0)=m\Delta c(0,m;A)=\frac{A_0(1/m)-1}{1/m}\geq
A_0'(0)\geq-1
\]
and
\[
\bigl(A^*\bigr)'(1)=m\Delta c(m-1,m;A)=\frac{1-A_0(1-1/m)}{1/m}\leq
A_0'(1)\leq1.
\]
This shows that $A^*\in\mathcal{A}$.
\end{pf}

The following proposition says that Bernstein approximations of Pickands
functions are themselves Pickands functions.
%

%
\begin{proposition}
\label{thm:BernsteinApproxPick}
We have $\mathcal{A}_m^+ \subset\mathcal{A}_m$, for every $m\geq1$.
Moreover, $\mathcal{A}^+_m=\mathcal{A}_m$ for $m=1,2,3$, while
$\mathcal{A}^+_4\neq\mathcal{A}_4$.
\end{proposition}

\begin{pf}
Let $V(t)=t\vee(1-t)$, $t\in[0,1]$. For the first statement, let
$A\in\mathcal{A}$ and $B_m(A,\cdot)\in\mathcal{A}_{m}^+$,
$m\geq1$. We have $B_m(A,0)=A(0)=1=A(1)=B_m(A,1)$, and $B_m(A,\cdot)$
is convex, see
Section~\ref{sect:BernsteinApprox}. Using Jensen's inequality, we obtain
\[
V(t)\leq A(t)=A\bigl(\expec_t(S_m/m)\bigr)\leq
\expec_t\bigl(A(S_m/m)\bigr)=B_m(A,t),\qquad  t
\in[0,1],
\]
and this shows that the \textit{boundary conditions} of a Pickands
function are
verified and so $B_m(A,\cdot)\in\mathcal{A}_m$.

For the second statement, recall that
$\mathcal{A}_1=\mathcal{A}_0=\{1\}=\mathcal{A}_1^+$. If
$A\in\mathcal{A}_m$ for either $m=2$ or $m=3$, then the arguments presented
in Example~\ref{ex:A2A3} can serve to verify that the piecewise
linear function $A^*$
interpolating the points $\{(k/m,c(k,m;A))\}_{k=0}^m$ belongs to
$\mathcal{A}$. Lemma~\ref{thm:PiecewisePick} then implies that
$A\in\mathcal{A}^+_m$, when $m=2,3$. The fact that $A^*$ belongs to
$\mathcal{A}$ is not necessarily true for $m=4$ as the following
shows: take
for instance
%
%
\begin{equation}
\label{eq:A4+} A(t)=1-t(1-t)\bigl\{1-2t(1-t)\bigr\},\qquad  t\in[0,1].
\end{equation}
It can be easily
verified that $A\in\mathcal{A}_4$, and we have
$A(t)=\expec_t\{c(S_4,4;A)\}$, with $c(k,4;A)=1$, for $k=0,2,4$ and
$c(1,4;A)=3/4=c(3,4;A)$. Suppose that
$A(\cdot)=B_4(A_0,\cdot)$, for some $A_0\in\mathcal{A}$. Then
$1=c(2,4;A)=c(2,4;B_4(A_0,\cdot))=A_0(1/2)$. This implies that
$A_0(t)=1$, for all $t\in[0,1]$, which in turn implies
that $1=B_4(A_0,t)=A(t)$, for all $t\in[0,1]$. Therefore,
$A\in\mathcal{A}_4\setminus\mathcal{A}_4^+$.
\end{pf}

\begin{Remark*} The above result says that
$\mathcal{A}^+_4\neq\mathcal{A}_4$, but we can say even more than this:
$\mathcal{A}^+_m\neq\mathcal{A}_m$ for every $m\geq4$. This follows
as a direct
application of \textit{The gap} theorem coming up (notice that for $A$
given by
in \eqref{eq:A4+} we have $A''(1/2)=0$).

For $m\geq0$, we call $\mathcal{A}_m$ the \textit{full model}
and $\mathcal{A}_m^+$ the \textit{submodel}. In order to give a
parameterization of the submodel, let
\[
\mathcal{H}_m^+ =\bigl\{A'': A\in
\mathcal{A}_{m+2}^+\bigr\} \subset\mathcal{H}_{m}.
\]
The next theorem shows that the (parametric) space of the coefficients
\[
C_m^+=\bigl\{ \bigl(c(0,m;h),\ldots,c(m,m;h)\bigr): h\in
\mathcal{H}^+_m\bigr\},
\]
is a polytope by showing
that $\mathcal{H}_m^+=\{h\in\mathcal{H}_m : c(k,m;h)\geq0,
k=0,\ldots,m\}$,
$m\geq0$.
\end{Remark*}
%

%
\begin{theorem}[(Parameterization)]
\label{thm:H+}
We have $h\in\mathcal{H}^+_m$, $m\geq0$, if and only if
%
%
\begin{equation}
\label{eq:boundary} \frac{1}{m+1} \sum_{j=0}^m
\biggl(1-\frac{j+1}{m+2} \biggr)c(j,m;h)\leq1, \qquad\frac
{1}{m+1} \sum
_{j=0}^m \frac{j+1}{m+2} c(j,m;h)\leq1
\end{equation}
and
%
%
\begin{equation}
\label{eq:positivity} c(k,m;h)\geq0,\qquad  k=0,\ldots,m.
\end{equation}
\end{theorem}

\begin{pf}
Let $m \geq0$.
If $h\in\mathcal{H}^+_m$ then $h\in\mathcal{H}_m$ and, by
definition of
$\mathcal{H}_m$,
the conditions
\[
\frac{1}{m+1} \sum_{j=0}^m \biggl(1-
\frac{j+1}{m+2} \biggr)c(j,m;h)\leq1, \qquad\frac{1}{m+1} \sum
_{j=0}^m \frac{j+1}{m+2} c(j,m;h)\leq1,
\]
hold. When $h\in\mathcal{H}^+_m$, there exists $A\in\mathcal
{A}_{m+2}$ such that
$h(\cdot)=B''_{m+2}(A,\cdot)$.
Therefore, by Lemma~\ref{thm:HtoA} and convexity of $A$,
\begin{eqnarray*}
c(k,m;h)&=&(m+2) (m+1) \biggl\{A \biggl(\frac{k+2}{m+2} \biggr)\\
&&\hphantom{(m+2) (m+1) \biggl\{}{}-2A
\biggl(
\frac{k+1}{m+2} \biggr)+A \biggl(\frac{ k } { m+2 } \biggr) \biggr
\} \geq0,\qquad
k=0,\ldots,m.
\end{eqnarray*}
Conversely, suppose that a polynomial $h\in\mathcal{P}_m$
satisfies \eqref{eq:boundary} and \eqref{eq:positivity}, then
$h\in\mathcal{H}_m$, and from Theorem~\ref{thm:representation},
there exists an $A\in\mathcal{A}_{m+2}$ such that $h=A''$. Now let
$A^*:[0,1]\to\mathbb{R}$ be the piecewise linear interpolation
of $\{(k/(m+2),c(k,m+2;A)\}_{k=0}^{m+2}$. The \textit{endpoint
conditions} on
the coefficients of $A$ given by Proposition~\ref{thm:pol_constraints} directly
imply that $A^*$ satisfies the \textit{endpoint conditions} of a Pickands
function. The convexity of $A^*$ is obtained by the
nonnegativity of the coefficients of $h$: Lemma~\ref{thm:HtoA} gives
$(m+2)(m+1)\Delta^2 c(k,m+2;A)=c(k,m;h)\geq0$, $k=0,\ldots,m$, and
the first
part of Lemma~\ref{thm:PiecewisePick} gives the convexity of $A^*$.
Therefore, $A^*\in\mathcal{A}$, and by the second part of
Theorem~\ref{thm:BernsteinApproxPick}, $A\in\mathcal{A}_{m+2}^+$,
which means
that $h\in\mathcal{H}_m^+$.
\end{pf}

A useful property of the submodel $\mathcal{A}_m^+$, $m\geq
1$, is that it is nested. We show this in the next theorem, and we find
the gap
between the full model and the submodel.
The proof relies on the Lorentz degree of a positive
polynomial on $(0,1)$, see for instance Powers and Reznick \cite{PR00}.
In the Bernstein basis, it is easy to see that $P\in\mathcal{P}_m$ is
nonnegative if $c(k,m;P)\geq0$, $k=0,\ldots,m$. This sufficient
condition is
not a necessary condition, see Karlin and Shapley \cite{KS53},
Szeg{\H{o}} \cite{Szego75} and P{\'o}lya and Szeg{\H{o}} \cite
{PS98}. It
may happen that for some $M>m$,
\[
\min\bigl\{c(k,m;P): k=0,\ldots,m\bigr\}<0, \qquad \mbox{while } \min
\bigl\{c(k,M;P):
k=0,\ldots,M\bigr\}\geq0.
\]
The question is when does this happen? For the (interesting) case:
$\deg(P)>0$, the necessary and sufficient condition for having $\min\{
c(k,M;P):
k=0,\ldots,M\}\geq0$ for some $M\geq\deg(P)$ is that
$P$ be positive on $(0,1)$. This is known in the literature as Bernstein's
theorem, and it motivates the definition of the Lorentz degree of a
polynomial $P\in\mathcal{P}_m$ which is positive on $(0,1)$:
\[
r(P)=\min\bigl\{M\geq m: c(k,M;P)\geq0, k=0,\ldots,M\bigr\}.
\]
To illustrate this notion, here is a very instructive example that will
lead us
naturally to \textit{The gap} theorem. The symmetric polynomials
in $\mathcal{A}_4\setminus\mathcal{A}_2$ have the following form:
\begin{eqnarray*}
&&A_{\alpha,\beta}(t)=1-\alpha t(1-t)\bigl\{1-\beta t(1-t)\bigr\},\\
&&\quad t\in[0,1],
\mbox{ with } \alpha\in(0,1] \mbox{ and } \beta\in[-1,2]\setminus
\{0\}.
\end{eqnarray*}
As can be seen in Figure~\ref{fig:Asym} below, $\alpha$ (which is the absolute
value of the first derivative at the endpoints) sets the triangle in
which the
curves belong, and then $\beta$ controls their curvature.
%
%
\begin{figure}

\includegraphics{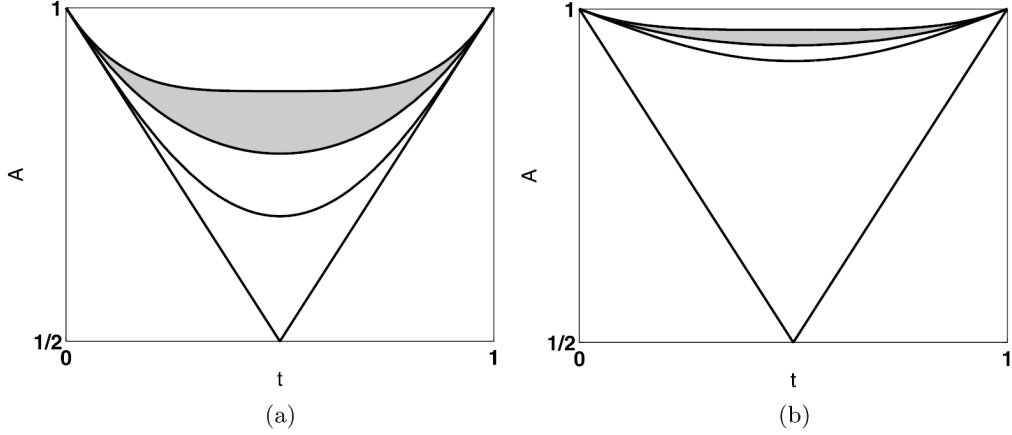}

\caption{Illustration of the gap between $\mathcal{A}_4^+$ and
$\mathcal{A}_4$
in the symmetric case. Here $A=A_{\alpha,\beta}$ with
$A_{\alpha,\beta}(t)=1-\alpha t(1-t)[1-\beta t(1-t)]$,
$t\in[0,1]$. In both plots, the three curves correspond to $\beta=-1,
1/2$ and 2
respectively (from bottom to top). The shaded regions hold the curves in
$\mathcal{A}_4\setminus
\mathcal{A}_4^+$. (a) $\alpha=1$. (b)~$\alpha=1/4$.}
\label{fig:Asym}
\end{figure}
Let $h_{\alpha,\beta}=A_{\alpha,\beta}''\in\mathcal{H}_2$, and
lets find the
values of $m\geq2$, such that $h_{\alpha,\beta}\in\mathcal
{H}_m^+$. The
smallest such $m$, if it exists, is the Lorentz degree of $h_{\alpha
,\beta}$,
and as we will see, it depends on $\beta$ here.
The function $h_{\alpha,\beta}$ is given, for $t\in[0,1]$, by
\begin{eqnarray*}
h_{\alpha,\beta}(t)&=&2\alpha\bigl\{(1+\beta)-6\beta t(1-t)\bigr\}\\
&=&2\alpha
\bigl[(1+\beta)\bigl\{(1-t)+t\bigr\}^m-6\beta t(1-t)\bigl\{(1-t)+t
\bigr\}^{m-2}\bigr],\qquad  m\geq2.
\end{eqnarray*}
Using Newton's binomial formula, we get
%
%
\begin{equation}
\label{all_coef} c(k,m;h_{\alpha,\beta})=2\alpha\biggl\{(1+\beta
)-6\beta
\frac
{k(m-k)}{m(m-1)} \biggr\},\qquad  k=0,\ldots,m, m\geq2.
\end{equation}
Fix $\alpha\in(0,1]$. From \eqref{all_coef}, we see that when $\beta
\in
[-1,0)$, we have $c(k,m;h_{\alpha,\beta})\geq0$, for $k=0,\ldots, m$,
and this
holds for every $m\geq2$. Therefore, the Lorentz degree in this case is
$r(h_{\alpha,\beta})=2$, and $h_{\alpha,\beta}\in\mathcal{H}_m^+$
if and only
if $m\geq2$. Next, when $\beta\in(0,2)$ the smallest coefficient
occurs when
$k=\lfloor m/2\rfloor$, and we get
%
%
\begin{eqnarray}
\label{min_coef} &&\min\bigl\{c(k,m;h_{\alpha,\beta}): k=0,\ldots
,m\bigr\}\nonumber\\[-8pt]\\[-8pt]
&&\quad =2\alpha
\biggl\{(1+\beta)-6\beta\frac{\lfloor
m/2\rfloor(m-\lfloor m/2\rfloor)}{m(m-1)} \biggr\},\qquad  m\geq2.\nonumber
\end{eqnarray}
By looking at the sign of \eqref{min_coef} and solving the inequality
for $m$
we obtain $c(\lfloor m/2 \rfloor,m;\break h_{\alpha,\beta}) \geq0$ if and
only if
$m\geq2\lceil(1+\beta)/(2-\beta)\rceil$. This tells us that the Lorentz
degree in this case is
\[
r(h_{\alpha,\beta})=2\bigl\lceil(1+\beta)/(2-\beta)\bigr\rceil,\qquad
\beta\in(0,2),
\]
and also that $h_{\alpha,\beta}\in\mathcal{H}_m^+$ if and only if
$m\geq r(h_{\alpha,\beta})$. The Lorentz degree increases without
bound as
$\beta$ increases to 2. Finally when $\beta=2$, we get
\[
c\bigl(\lfloor m/2\rfloor,m;h_{\alpha,2}\bigr)=\frac{-6\alpha
}{2\lceil
m/2\rceil
-1}<0, \qquad \mbox{for every } m\geq2,
\]
so $h_{\alpha,2}\in\mathcal{H}_2\setminus\bigcup_{m\geq2}\mathcal
{H}_m^+$ or
equivalently
$A_{\alpha,2}\in\mathcal{A}_4\setminus\bigcup_{m\geq4}\mathcal
{A}_m^+$. As the
following theorem shows, this is true because $h_{\alpha,2}(1/2)=0$,
that is,
$h_{\alpha,2}$ is nonnegative but fails to be positive on $(0,1)$.
Finally, for this example, we note that $A_{\alpha,\beta}\in\mathcal{A}_4^+$
if and only if $(\alpha,\beta)\in(0,1]\times[-1,0)\cup(0,1/2]$, see
Figure~\ref{fig:Asym}.
%

%
\begin{theorem}[(The gap)]
\label{thm:Gap}
For all $m\geq1$, we have $\mathcal{A}_m^+ \subset\mathcal
{A}_{m+1}^+$. Moreover,
for $h\in\bigcup_{m=0}^{\infty}\mathcal{H}_m$ we have $h\notin\bigcup
_{m=0}^{\infty}\mathcal{H}^+_m$
if and only if $\deg(h)>0$ and $h(t)=0$ for some $t\in(0,1)$.
\end{theorem}

\begin{pf}
First, since $\mathcal{A}_m^+=\mathcal{A}_m$, for $m=1,2,3$, we
already have
$\mathcal{A}^+_1\subset\mathcal{A}^+_2\subset\mathcal{A}^+_3$.
Let $h\in\mathcal{H}_{m}^+$, for some $m\geq1$.
We want to show that $h\in\mathcal{H}_{m+1}^+$.
Since $\mathcal{H}_m^+\subset\mathcal{H}_{m}\subset\mathcal
{H}_{m+1}$ we just
need to show that
$c(k,m+1;h)\geq0$, $k=0,\ldots,m+1$.
Using the binomial identities from Section~\ref{sect:binident}, we have
\begin{eqnarray*}
&& \expec_t\bigl\{c(S_{m+1},m+1;h)\bigr\} 
\\
&&\quad =
\expec_t\bigl\{c(S_m,m;h)\bigr\}
\\
&&\quad = \expec_t\bigl\{(1-t)c(S_m,m;h)+tc(S_m,m;h)
\bigr\}
\\
&&\quad = \frac{1}{m+1} \expec_t\bigl\{
(m+1-S_{m+1})c(S_{m+1},m;h)+S_{m+1}c(S_{m+1}-1,m;h)
\bigr\},
\end{eqnarray*}
leading to the identity
\[
c(j,m+1;h)=\frac{j}{m+1}c(j-1,m;h)+ \biggl(1-\frac{j}{m+1}
\biggr)c(j,m;h),\qquad  j=0,\ldots,m+1.
\]
From Theorem~\ref{thm:H+} it follows that $h\in\mathcal{H}^+_{m+1}$, which
proves the first assertion.

Now let $h\in\bigcup_{m=0}^{\infty}\mathcal{H}_m$. To prove the second
assertion, we
show that we have $h\in\bigcup_{m=0}^{\infty}\mathcal{H}^+_m$ if and
only if
$\deg(h)=0$ or $h(t)>0$ for all $t\in(0,1)$.

First, assume that $h\in\bigcup_{m=0}^{\infty}\mathcal{H}^+_m$. If
$h\in\mathcal{H}^+_0$,
then $\deg(h)=0$. Otherwise, $h\in\mathcal{H}^+_m\setminus\mathcal
{H}^+_0$, so
$\deg(h)>0$ and this implies that $\prob_t\{c(S_m,m;h)\geq0\}=1$ and
$\prob_t\{c(S_m,m;h)> 0\}>0$ for all $t\in(0,1)$. It follows that
$h(t)=\expec_t\{c(S_m,m;h)\}>0$, for all $t\in(0,1)$.

Conversely, if $\deg(h)=0$ then $h\in\mathcal{H}^+_0$, while if $h$
is positive on
$(0,1)$ then $h\in\mathcal{H}^+_{r(h)}$, where $r(h)$ is the Lorentz
degree of $h$.
\end{pf}

\section{On the quality of the Bernstein approximations}
\label{sect:quality}
We are now concerned with the flexibility of the submodel $\mathcal{A}_m^+$
consisting of the Bernstein approximations of Pickands functions. More
precisely, we provide an answer to how well a Pickands function
$A\in\mathcal{A}$ can be approached in the space $\mathcal{A}_m^+$.
We also
work out the range $\tau(\mathcal{A}_m^+)$, as $m$ varies, of some dependence
measures $\tau$ on $\mathcal{A}$.
%

%
\begin{theorem}[(Approximation)]
\label{thm:BernsteinApprox}
For $m\geq1$, $A\in\mathcal{A}$, $t\in[0,1]$, we have
%
%
\begin{equation}
\label{eq:BernapproxA} A(t)\leq B_m(A, t)\leq A(t)+2t(1-t)
\prob_t\bigl(S_{m-1}=\lfloor mt \rfloor\bigr),
\end{equation}
where
\[
2t(1-t)\prob_t\bigl(S_{m-1}=\lfloor mt \rfloor\bigr)=
\biggl\{\frac
{2t(1-t)}{m\pi} \biggr\}^{1/2}+\mathrm{O}\bigl(m^{-3/2}
\bigr),\qquad  (m\to\infty).
\]

Moreover, when $A=V$, with $V(t)=(1-t)\vee t$, $t\in[0,1]$, we get the finer
approximation
%
%
\begin{equation}
\label{eq:BernapproxV} B_m(V, t)-V(t) \leq\bigl\{1-V(t)\bigr\}
\prob_t\bigl(S_{m-1}=\lfloor m/2 \rfloor\bigr),
\end{equation}
where $\{1-V(t)\}\prob_t(S_{m-1}=\lfloor m/2\rfloor)=0$, for $t\in\{
0,1\}$,
while
\[
\bigl\{1-V(t)\bigr\}\prob_t\bigl(S_{m-1}=\lfloor m/2
\rfloor\bigr)\leq\frac{1}{2} \biggl\{\frac{1}{t}\wedge
\frac
{1}{(1-t)} \biggr\} \biggl\{ \frac{ 2t(1-t)}{m\pi} \biggr\}^{1/2}+
\mathrm{O}\bigl(m^{-3/2}\bigr), \qquad (m\to\infty),
\]
for $t\in(0,1)$, with equalities in the last two inequalities when $t=1/2$.
\end{theorem}

\begin{pf}
The inequality $B_m(A, \cdot) \geq A$ follows from Jensen's inequality
\[
B_m(A, t)=\expec_t\bigl\{A(S_{m}/m)\bigr\}
\geq A(t),\qquad  t\in[0,1].
\]

The approximation formula \eqref{eq:BernapproxA} is obtained as follows.
First, the convexity and the fact that $V(t)\leq A(t)\leq1$, $t\in[0,1]$
imply that
$|A(t_1)-A(t_2)|\leq|t_1-t_2|$, for all $t_1,t_2 \in[0,1]$. From
that, for
$t\in[0,1]$,
\begin{eqnarray*}
B_m(A, t)-A(t)&=&\bigl|\expec_t\bigl\{A(S_{m}/m)-A(t)
\bigr\}\bigr|\\
&\leq& \expec_t\bigl\{\bigl|A(S_{m}/m)-A(t)\bigr|\bigr\}
\\
&\leq& \expec_t\bigl\{|S_{m}/m-t|\bigr\}
\\
&=& 2t(1-t)\prob_t\bigl(S_{m-1}=\lfloor mt \rfloor\bigr)
\\
&=& \biggl\{\frac{2t(1-t)}{m\uppi} \biggr\}^{1/2}+\mathrm{O}
\bigl(m^{-3/2}\bigr), \qquad (m\to\infty),
\end{eqnarray*}
the last two equalities can be found in, for instance, Johnson \cite{Johnson57}. In
particular, the last equality is obtained using Stirling's formula,
$n!=(2\uppi
n)^{1/2}(n/\mathrm{e})^n\{1+\mathrm{O}(1/n)\}$, ($n\to\infty$). Here we
should point
out that the
Cauchy
inequality gives the weaker result
\[
B_m(A, t)-A(t)\leq\expec_t\bigl(|S_m/m-t|\bigr)\leq
\bigl\{t(1-t)/m\bigr\}^{1/2},
\]
although the rate of convergence remains of the same order.

To show \eqref{eq:BernapproxV}, making use of the binomial identities (once
again), we have,
for
$t\in[0,1]$,
\begin{eqnarray*}
B_m(V, t)&=&\expec_t\bigl\{V(S_{m}/m)\bigr\}\\
&=&
\expec_t\bigl\{(1-S_{m}/m)\indic(S_{m}\leq m/2)
+(S_{m}/m)\indic(S_{m}> m/2)\bigr\}
\\
&=&(1-t)\prob_t(S_{m-1}\leq m/2)+t\prob_t(S_{m-1}>
m/2-1).
\end{eqnarray*}
Thus,
\begin{eqnarray*}
B_m(V, t)-V(t)&=&\bigl\{1-V(t)\bigr\}\prob_t(m/2-1<S_{m-1}
\leq m/2)
\\
&&{}+\bigl\{1-t-V(t)\bigr\}\prob_t(S_{m-1}\leq m/2-1)+ \bigl
\{t-V(t)\bigr\}\prob_t(S_{m-1}>m/2)
\\
&=&\bigl\{1-V(t)\bigr\}\prob_t(m/2-1<S_{m-1}\leq m/2)
\\
&&{}-|2t-1|\bigl\{\prob_t(S_{m-1}\leq m/2-1)\indic(t> 1/2)+
\prob_t(S_{m-1}>m/2)\indic(t\leq1/2)\bigr\}
\\
&\leq&\bigl\{1-V(t)\bigr\}\prob_t(m/2-1<S_{m-1}\leq m/2)
\\
&=&\bigl\{1-V(t)\bigr\}\prob_t\bigl(S_{m-1}=\lfloor m/2
\rfloor\bigr).
\end{eqnarray*}
The latter inequality being an equality at $t=1/2$.

Finally, $\prob_t(S_{m-1}=\lfloor m/2 \rfloor)\leq\prob
_t(S_{m-1}=\lfloor mt
\rfloor)$,
and the rate of convergence follows again by Stirling's formula.
\end{pf}

Dependence measures for bivariate extremes have been studied in the
literature, see for instance Tawn \cite{Tawn88} and Weissman \cite{Weissman08}. In
particular, for $A\in\mathcal{A}$ the following measures were proposed
\[
\tau_1(A)=2\bigl\{1-A(1/2)\bigr\},\qquad  \tau_2(A)=\expec
\bigl[4\bigl\{1-A(U)\bigr\}\bigr], \qquad U\sim\mathcal{U}(0,1).
\]
Let $A\in\mathcal{A}$ and
consider the regions $R_i\subset\mathbb{R}^2$, $i=1,2,3$, given by
\begin{eqnarray*}
R_1 &=& \bigl\{(t,y)\colon0\leq t\leq1, A(1/2)+\bigl(1-A(1/2)\bigr)
\mid2t-1\mid\leq y \leq1\bigr\},
\\
R_2 &=& \bigl\{(t,y)\colon0\leq t\leq1, A(t)\leq y \leq1\bigr\},
\\
R_3 &=& \bigl\{(t,y)\colon0\leq t\leq1, V(t)\leq y \leq1\bigr\}.
\end{eqnarray*}
Since $A\in\mathcal{A}$ we have $R_1\subset R_2\subset R_3$ and
$0\leq
\tau_1\leq\tau_2\leq1$ because,
\[
\tau_i=4\times\area(R_i),\qquad  i=1,2 \quad \mbox{and}\quad  0\leq
\area(R_1)\leq\area(R_2)\leq\area(R_3)=1/4.
\]
In particular,
\[
\tau_1\bigl\{B_m(A, \cdot)\bigr\}=\expec_{1/2}
\bigl[2\bigl\{1-A(S_m/m)\bigr\}\bigr],\qquad  \tau_2\bigl
\{B_m(A, \cdot)\bigr\}=\expec\bigl[4\bigl\{1-A(U_m/m)\bigr
\}\bigr],
\]
where $U_m\sim\mathcal{U}\{0,\ldots,m\}$. Note that since $S_{m}/m$ and
$U_m/m$ converge in distribution to $t$ and $U\sim\mathcal{U}(0,1)$
respectively, and since $A$ is continuous and bounded, we obtain that
$\tau_i\{B_m(A, \cdot)\}\to\tau_i(A)$, as $m\to\infty$, $i \in\{
1,2\}$.
The next proposition gives the
range of these measures for the Bernstein approximations.
%

%
\begin{proposition}[(Dependence measures)]
\label{thm:depMeasures}
Let $V(t)=(1-t)\vee t$, $t\in[0,1]$ (the Pickands
function corresponding to complete monotone dependence). For $m\geq1$,
$\tau_i(\mathcal{A}_{m}^+)=[0,\tau_i\{B_m(V, \cdot)\}]$, $i \in
\{1,2\}$, with
\[
\tau_1\bigl\{B_m(V, \cdot)\bigr\}=1-\prob_t
\bigl(S_{m-1}=\lfloor m/2 \rfloor\bigr),\qquad  \tau_2\bigl
\{B_m(V, \cdot)\bigr\}=\lfloor m/2\rfloor/ \bigl(\lfloor m/2\rfloor+
1/2\bigr).
\]
\end{proposition}

\begin{pf}
First, if $A_1\leq A_2$, $A_1,A_2\in\mathcal{A}$, then $B_m(A_1,\cdot
)\leq
B_m(A_2, \cdot)$, so that $\tau_i\{B_m(A_1, \cdot)\}\geq\tau_i\{B_m(A_2,
\cdot)\}$, $i=1,2$, and
this implies that $\tau_i(\mathcal{A}_{m}^+)\subset[0,\tau_i\{
B_m(V, \cdot)\}]$.
The reverse inclusion follows from the convexity of the space $\mathcal
{A}_{m}^+$,
the fact that $B_m(V, \cdot)$ is a Pickands function and the
linearity (under convex combinations) of the functionals $\tau_1$ and
$\tau_2$.\looseness=1
\end{pf}

\section{Simulation experiment}
\label{sect:simulations}
We now compare the maximum likelihood estimator from the full
model and the one from the submodel through simulated data. In both
cases, the
maximum likelihood estimator (MLE) is computed numerically. For a fixed
$m\geq
0$, when using the full model, the estimator is the polynomial
$A_{\hat{\theta}}$, where $\hat{\theta}$ is the
solution of the nonlinear constrained maximization problem
%
%
\begin{equation}
\label{eq:mleF} \hat{\theta}=\mathop{\argmax}_{\theta\in\Theta_m}
L_F(\theta\mid
u,v),
\end{equation}
where $\Theta_m$ is given in \eqref{eq:Thetam}, for $m\geq1$,
$\Theta_0=[0,2]$,
$(u,v)=\{(u_1,v_1),\ldots,(u_n,v_n)\}$ is
the data and the likelihood $L_F$ is obtained by computing the coefficients
of $h$ via Lemma~\ref{thm:PandQ}, by deriving its associated Pickands function
$A$ via Lemma~\ref{thm:HtoA}, and then by using the mixed partial derivative
of \eqref{eq:CA}. The estimator for the Bernstein approximations submodel
is $A_{\hat{h}}$, the problem to solve is
%
%
\begin{equation}
\label{eq:mleS} \hat{c}=\mathop{\argmax}_{c\in C_m^+} L_S(c \mid u,v),
\end{equation}
where $C_m^+$ is
the polytope given in
Theorem~\ref{thm:H+} and $\hat{h}\in\mathcal{H}_m^+$ is such that
$\hat{c}=(c(k,m;\hat{h}))_{k=0}^m$.
Again the likelihood $L_S$ is obtained by Lemma~\ref{thm:HtoA}, and
\eqref{eq:CA}.
Both problems have solutions by continuity of the likelihood over
their (compact) feasible regions $\Theta_m$ and $C_m^+$, respectively. These are easily solved by using the \matlab global
search algorithm, for instance.

We consider the three following models for the simulation: the asymmetric
logistic model
%
%
\begin{equation}
\label{eq:alog} A(t)=(1-\psi_1)t+(1-\psi_2) (1-t)+\bigl[(
\psi_1t)^{1/\alpha}+\bigl\{\psi_2(1-t)\bigr
\}^{
1/\alpha}\bigr]^{\alpha}, \qquad t\in[0,1],
\end{equation}
with $\alpha\in(0,1]$, $0\leq\psi_1,\psi_2 \leq1$,
the symmetric mixed model
%
%
\begin{equation}
\label{eq:mix} A(t)=1-\psi t+ \psi t^2, \qquad t\in[0,1],
\end{equation}
$\psi\in[0,1]$, and the polynomial $A$ obtained via Lemma~\ref
{thm:HtoA} with
$c(0,2;h)= 2$, $c(1,2;h)=-1/3$ and $c(2,2;h)=1/5$,
leading to
%
%
\begin{equation}
\label{eq:polfull} A(t)= 1- (83/180)t+ t^2- (7/9) t^3+(43/180)t^4,
\qquad t\in[0,1].
\end{equation}
While it is clear from the characterization of
$\mathcal{H}_2$ given after the proof of
Lemma~\ref{thm:HtoA} that $A\in\mathcal{A}_{4}\setminus\mathcal
{A}_{4}^+$, it turns out that it has Lorentz
degree $r(h)=6$, where $h=A''$, so that $A \in\mathcal
{A}_{8}^+\setminus\mathcal{A}_{7}^+$. The first two models have also
been considered for simulations
studies in Fils-Villetard \textit{et al.} \cite{FGS08} and B{\"u}cher \textit{et al.}
\cite{BDV11}. Here we take the parameter values
$\alpha=1/2$, $\psi_1=1/10$, $\psi_2=1/2$ in the first model, and
$\psi=9/10$ in
the second model. Note that for this particular choice of asymmetric logistic
model, $A\notin\bigcup_{m\geq0} \mathcal{A}_{m}$, and for the
symmetric mixed model,
$A \in\mathcal{A}_{2}^+=\mathcal{A}_{2}$.

\begin{figure}[t]

\includegraphics{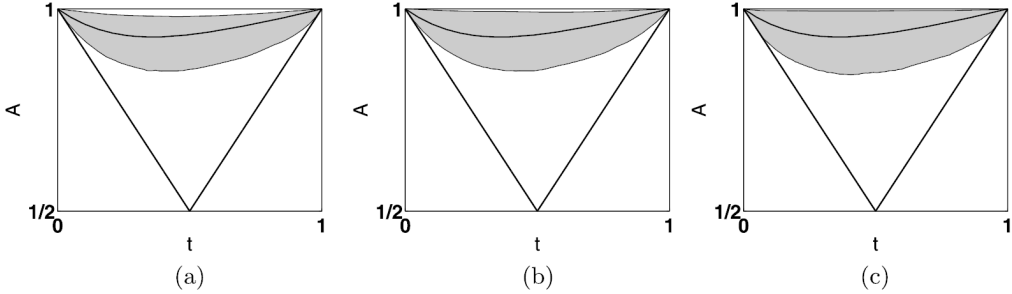}

\caption{Shaded region corresponds to point-wise 95\% confidence
intervals for
the true Pickands function (black curve) of model \protect\eqref{eq:polfull}.
On part (a) are the results using the maximum
likelihood estimator from the full model computed via \protect\eqref
{eq:mleF}, in
part (b) using the maximum likelihood estimator from the
submodel computed via \protect\eqref{eq:mleS} and on part (c)
using the optimal
CFG estimator \protect\eqref{eq:CFGopt}. Here, $m=5$ for the polynomials. (a) Model~(\ref{eq:polfull}), full MLE $A_{\hat{\theta}}$.
(b) Model~(\ref{eq:polfull}), sub MLE $A_{\hat{h}}$. (c) Model~(\ref{eq:CFGopt}), $\hat{A}_{\CFGopt}$.}
\label{fig:band}
\end{figure}

\begin{figure}

\includegraphics{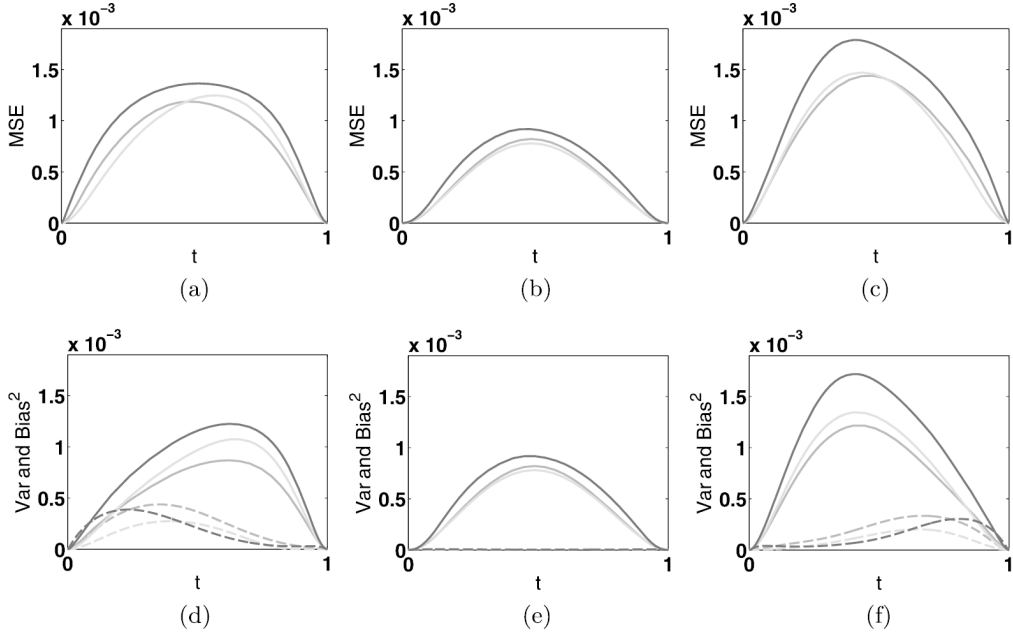}

\caption{In all the above illustrations, different shades of grey represent
(from dark to light grey): the optimal CFG estimator, the MLE
from the full model, and the MLE from the submodel respectively. In parts (a)--(c), the curves represent the estimated mean squared
error of
the estimators. Parts (d)--(f) show the variances (thick curves) and squared bias
(dashed curves). Here, $m=5$ for the polynomials. (a) Model~(\ref{eq:alog}). (b) Model~(\ref{eq:mix}).
(c) Model~(\ref{eq:polfull}). (d) Model~(\ref{eq:alog}). (e) Model~(\ref{eq:mix}).
(f) Model~(\ref{eq:polfull}).}
\label{fig:MSE}
\end{figure}

A general algorithm for drawing independent couples from the copulas
associated to these models is provided by Ghoudi \textit{et al.} \cite
{GKR98}. Here we draw
$1 000$ samples of sizes $n=100$ from the copulas corresponding to each of
the above three models. For every sample, the maximum likelihood
estimate is
obtained using the full model by solving \eqref{eq:mleF} and the submodel
by solving \eqref{eq:mleS} both with $m=5$. This gives polynomial estimates
of degree at most 7, that is $A_{\hat{\theta}} \in\mathcal{A}_{7}$
for the full
model and $A_{\hat{h}} \in\mathcal{A}_{7}^+$ for the submodel. For
comparisons,
we also consider
the popular $A_{\CFG}$ estimator from Cap{\'e}ra{\`a} \textit{et al.} \cite
{CFG97}. The latter can
lead to estimates which are not genuine Pickands functions because they
do not satisfy either the \textit{boundary conditions} or the \textit
{convexity
condition}. A modification which leads to Pickands functions as
estimates is
$\hat{A}_{\CFG}= \mbox{greatest convex minorant of } 1\wedge
\{A_{\CFG}\vee V\}$,
where $V(t)=(1-t)\vee t$, $t \in[0,1]$. Another good estimator,
according to Genest and Segers \cite{GS2009}, is the optimally
corrected CFG estimator
$A_{\CFGopt}$ proposed initially by Segers \cite
{Segers08}. The estimator that
we use for the comparisons is
%
%
\begin{equation}
\label{eq:CFGopt} \hat{A}_{\CFGopt}= \mbox{greatest convex minorant
of } 1
\wedge\{A_{\CFGopt}\vee V\}.
\end{equation}
We show 95\% point-wise confidence intervals for model \eqref{eq:polfull}
in Figure~\ref{fig:band}.
%
%
To compare the performance of the various estimators in the cases considered
here, we look at estimates\vadjust{\goodbreak} of the mean squared error
\[
\mse_A\bigl\{\hat{A}(t)\bigr\}=\expec\bigl\{\hat{A}(t)-A(t)
\bigr
\}^2,\qquad  t\in[0,1],
\]
where $\hat{A}$ is an estimator and $A$ is the true
Pickands function. We also look at the variance and bias
separately. These results are plotted in Figure~\ref{fig:MSE}.
%
%
In the simulation results, the MLE of the submodel (uniformly)
outperforms the optimal CFG estimator \eqref{eq:CFGopt} in terms of mean
squared error and also variance. The MLE of the full model has uniformly
smaller variance than \eqref{eq:CFGopt}. The bias of all three
estimators is much smaller in the symmetric case than in the other two cases.
When $m=5$, the results do not clearly indicate an overall winner
between the
two maximum likelihood estimators for sample sizes $n=100$. This seems to
change, however (according to further simulations we have done), when
$m$ is
greater with respect to the data size $n$. This will be noticed in the
following.

\begin{figure}

\includegraphics{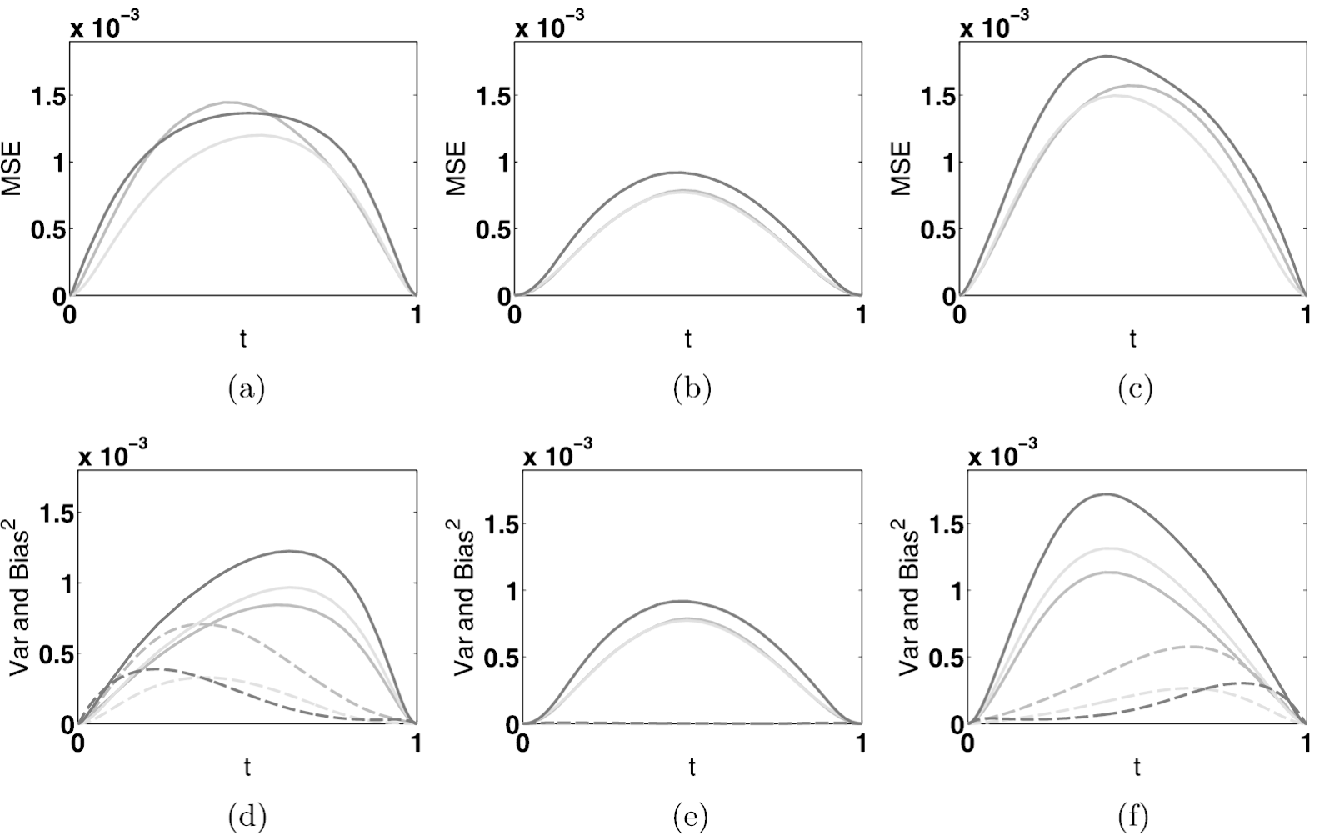}

\caption{In all the above illustrations, different shades of grey represent
(from dark to light grey): the optimal CFG estimator, the MLE
from the full model, and the MLE from the submodel respectively. In parts (a)--(c), the curves represent the estimated mean squared
error of
the estimators. Parts (d)--(f) show the variances (thick curves) and squared bias
(dashed curves). Here, $m=8$ for the polynomials. (a) Model~(\ref{eq:alog}).
(b) Model~(\ref{eq:mix}). (c) Model~(\ref{eq:polfull}). (d) Model~(\ref{eq:alog}).
(e) Model~(\ref{eq:mix}). (f) Model~(\ref{eq:polfull}).}
\label{fig:MSE8}
\end{figure}

We run the entire simulation again, using the same data set, but
this time with $m=8$. This gives polynomial estimates
of degree at most 10, that is $A_{\hat{\theta}} \in\mathcal
{A}_{10}$ for the
full model and $A_{\hat{h}} \in\mathcal{A}_{10}^+$ for the submodel.
The results
are plotted in Figure~\ref{fig:MSE8}.
%
%
We can observe that the
performance of the estimator from the submodel is roughly the same as
in the
case $m=5$, but the performance of the estimator from the full
model has somehow worsened in the two nonsymmetric cases in part due
to an
increase of bias.

\begin{figure}

\includegraphics{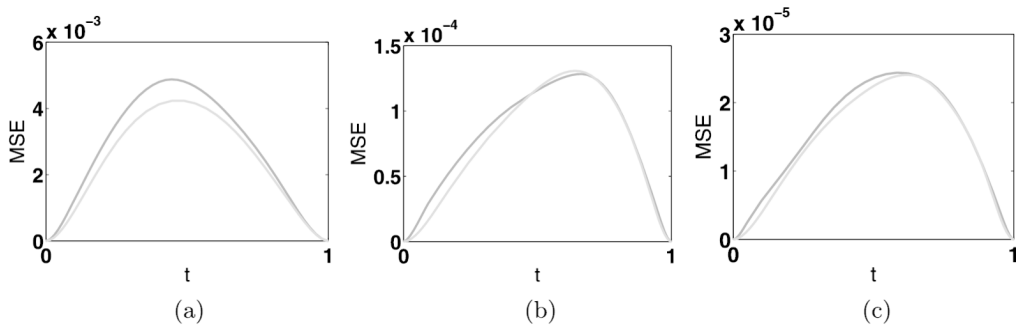}

\caption{Dark grey curve is the MSE of the MLE from the full model and
light grey curve is the MSE of the MLE from the submodel. Here, the true
Pickands function is \protect\eqref{eq:alog} and $m=5$ for the polynomials.
(a) $n=30$. (b) $n=1000$. (c) $n=5000$.}
\label{fig:MSEasymp}
\end{figure}

As a final numerical investigation, we
compare the performance of the MLE from the full model and the one
from the submodel, this time, when the sample size $n$ increases and
$m=5$ is
held fixed. Here, in every case we
have tried, the mean squared error of the MLE of the submodel was uniformly
smaller than that of the full model when the polynomial degree $m$ was
large in
comparison to the sample size $n$. Then, when $n$ increases and $m$ is held
fixed, the two estimators' mean squared error curves are close and cross
eachother. Figure~\ref{fig:MSEasymp} gives a typical illustration.
%
%

\section{Concluding remarks and comments}
\label{sect:conclusion}
\subsection*{The choice of basis}
First, the results obtained in
Section~\ref{sect:characterization} could be developed in any other polynomial
basis, and we took the popularity of the power basis into consideration.
However, it appears that some of the results have an easy
interpretation when the coefficients are expressed in the Bernstein basis
but seem meaningless when the coefficients are
expressed in the power basis. For instance, the
link between $A\in\mathcal{A}_{m+2}$ and
$h=A''$, together with the \textit{endpoint derivatives conditions}
can be
written
\begin{eqnarray*}
\label{eq:expec_representation} A(t)=1-\expec\bigl(\bigl[\bigl\{
(1-t)U\bigr\}\wedge\bigl\{t(1-U)
\bigr\}\bigr]h(U)\bigr),
\end{eqnarray*}
with $t\in[0,1]$, and
\[
\expec\bigl\{(1-U)h(U)\bigr\}\leq1, \qquad \expec\bigl\{Uh(U)\bigr\}
\leq1,
\]
with $U\sim\mathcal{U}(0,1)$. When we express these expressions in
terms of the
coefficients in the Bernstein basis, see
Lemma~\ref{thm:HtoA} and inequalities \eqref{eq:boundcondh}
respectively, we
get a striking
similarity
\begin{eqnarray*}
c(k,m+2;A)=1-\expec\bigl( \bigl[ \bigl\{ (1-t )U_m \bigr\} \wedge
\bigl\{t (1-U_m ) \bigr\} \bigr]c\bigl((m+2)U_m-1,m;h
\bigr) \bigr),
\end{eqnarray*}
with $t=k/(m+2)$, $k=0,\ldots, m+2$, and
\[
\expec\bigl\{(1-U_m)c\bigl((m+2)U_m-1,m;h\bigr)\bigr\}
\leq1, \qquad \expec\bigl\{U_mc\bigl((m+2)U_m-1,m;h\bigr)\bigr\}
\leq1,
\]
with $U_m\sim\mathcal{U}\{1/(m+2),\ldots,(m+1)/(m+2)\}$. Now let
\[
A(t)=\sum_{k=0}^{m+2} \alpha_k
t^k, \qquad h(t)=\sum_{k=0}^m
\eta_k t^k, \qquad t\in[0,1].
\]
If we do the same exercise in the power basis this time, we obtain
\begin{eqnarray*}
\alpha_0&=&1,\qquad  \alpha_1=-\sum
_{k=0}^m\frac{1}{(k+1)(k+2)}\eta_k,\\
\alpha_k&=&\frac{1}{k(k-1)}\eta_{k-2}, \qquad 2\leq k \leq m+2
\end{eqnarray*}
and
\[
\sum_{k=0}^m \frac{1}{(k+1)(k+2)}
\eta_k\leq1, \qquad \sum_{k=0}^m
\frac{1}{(k+2)}\eta_k\leq1,
\]
these are definitely difficult to grasp.

Secondly, in view of the model constructed in Section~\ref
{sect:approximation}, the Bernstein basis is clearly the right
basis for
approximating continuous functions on $[0,1]$ by polynomials, think
of the Weierstrass theorem for instance. It turns out that the
Bernstein basis
is very appealing for other reasons as well. For example, the geometry
of the
Pickands functions is directly reflected in the coefficients, this is
made clear
by Lemma~\ref{thm:PiecewisePick}. It is also shown that
$\mathcal{A}_m^+=\mathcal{A}_m$ for $m=1,2,3$, but
$\mathcal{A}_4^+\neq\mathcal{A}_4$. Working with the Bernstein
basis throughout
made this finding easier.

\subsection*{Full model vs submodel}
 The
convexity condition on the full model in Section~\ref
{sect:characterization} is
obtained via intermediate polynomials ($P$ and $Q$), so that the
coefficients of the resulting Pickands function $A$ are parameterized by
those of $P$ and $Q$. In this parameterization, there is an identifiability
problem with the parameters. To see this, in \eqref{eq:h}, take for
example the
four couples $(P,Q)$, $(-P,Q)$,
$(Q,-P)$ and
$(-P,-Q)$, these will all produce the same function $h$. However, for
inferential
purposes, the parameters of interest remain the coefficients of $h$ (or $A$),
not the
ones of $P$ and $Q$. The submodel is not concerned with this issue at
all, the parameters of the submodel
are the coefficients. The
simplicity of this model, its polytopal parameter space, the gap
theorem and
its flexibility
makes it more appealing, for us at least, than the full model in
Section~\ref{sect:characterization}. A practical consequence of the simplicity
of the
submodel is in finding the maximum likelihood estimator;
it is more than twice as fast, numerically, than for the complete model.
Finally, we mention that although the likelihood is not concave (and rather
complicated), local maxima could indeed make difficult the search for a global
maximum. However, using a
state-of-the-art \matlab global search algorithm, we have not
encountered problematic situations.

%
\begin{appendix}
\section*{Appendix}\label{app}
\begin{pf*}{Proof of Lemma \protect\ref{thm:PandQ}}
Consider $X_1,X_2,\ldots\,$, a sequence of independent $\Bernoulli(t)$
random variables.
Let $S_1,S_2,\ldots\,$, be the sequence of the cumulative sums, $
S_n=\sum_{k=1}^n X_k$.
Suppose that\vspace*{2pt} $m$ is even, $m>0$, $P\in\mathcal{P}_{m/2}$,
$Q\in\mathcal{P}_{(m-2)/2}$.
We can write, by equation \eqref{eq:h},
\[
h(t)=P^2(t)+t(1-t)Q^2(t),\qquad  t\in[0,1].
\]
Now,
\begin{eqnarray*}
P^2(t) &=& \expec_t \bigl\{p(S_{m/2}
,m/2)p(S_m-S_{m/2}, m/2) \bigr\}
\\
&=& \expec_t\bigl[\expec\bigl\{p(Y,m/2)p(S_m-Y, m/2)
\mid S_m \bigr\}\bigr],
\end{eqnarray*}
with $\mathcal{L}(Y\mid S_m=k)=\Hypergeo(k,m/2,m)$, $k=0,1,\ldots,m$.
Similarly,
%
%
\begin{eqnarray}
\label{eq:Q2} Q^2(t) &=& \expec_t \bigl\{q
\bigl(S_{(m-2)/2} ,(m-2)/2\bigr)q\bigl(S_{m-2}-S_{(m-2)/2},
(m-2)/2\bigr) \bigr\}
\nonumber
\\[-8pt]\\[-8pt]
&=& \expec_t\bigl[\expec\bigl\{q\bigl(Y ,(m-2)/2\bigr)q
\bigl(S_{m-2}-Y, (m-2)/2\bigr)\mid S_{m-2} \bigr\}\bigr],\nonumber
\end{eqnarray}
with $\mathcal{L}(Y\mid S_{m-2}=k)=\Hypergeo(k,(m-2)/2,m-2)$,
$k=0,1,\ldots,m-2$.
From \eqref{eq:Q2} and the binomial identities, we get
\[
t(1-t)Q^2(t)= \expec_t \biggl[\frac{S_m(m-S_m)}{m(m-1)} \expec
\bigl\{q\bigl(Y_2,(m-2)/2\bigr)q\bigl(S_{m}-Y_2-1,
(m-2)/2\bigr)\mid S_m \bigr\} \biggr],
\]
with $\mathcal{L}(Y_2\mid S_{m}=k)=\Hypergeo(k-1,(m-2)/2,m-2)$,
$k=1,\ldots,m-1$.

Now, when $m$ is odd,
$P,Q\in\mathcal{P}_{(m-1)/2}$, and we can write
\[
h(t)=tP^2(t)+(1-t)Q^2(t),\qquad  t\in[0,1].
\]
Here,
%
%
\begin{eqnarray}
\label{eq:P2} P^2(t) &=& \expec_t \bigl\{p
\bigl(S_{(m-2)/2} ,(m-1)/2\bigr)p\bigl(S_{m-1}-S_{(m-1)/2},
(m-1)/2\bigr) \bigr\}
\nonumber
\\[-8pt]\\[-8pt]
&=& \expec_t \bigl[\expec\bigl\{p\bigl(Y,(m-1)/2\bigr)p
\bigl(S_{m-1}-Y, (m-1)/2\bigr)\mid S_{m-1} \bigr\} \bigr],\nonumber
\end{eqnarray}
with $\mathcal{L}(Y\mid S_{m-1}=k)=\Hypergeo(k,(m-1)/2,m-1)$,
$k=0,1,\ldots,m-1$.
Moreover, from \eqref{eq:P2} and the binomial identities,
\[
tP^2(t)= \expec_t \biggl[\frac{S_m}{m}\expec\bigl
\{p\bigl(Y_1,(m-1)/2\bigr)p\bigl(S_{m}-Y_1-1,
(m-1)/2\bigr)\mid S_{m} \bigr\} \biggr],
\]
with $\mathcal{L}(Y_1\mid S_{m}=k)=\Hypergeo(k-1,(m-1)/2,m-1)$,
$k=1,\ldots,m$.

Similarly,
\begin{eqnarray*}
Q^2(t) &=& \expec_t \bigl[\expec\bigl\{q\bigl(Y,(m-1)/2
\bigr)q\bigl(S_{m-1}-Y, (m-1)/2\bigr)\mid S_{m-1} \bigr\}
\bigr],
\end{eqnarray*}
with $\mathcal{L}(Y\mid S_{m-1}=k)=\Hypergeo(k,(m-1)/2,m-1)$,
$k=0,1,\ldots,m-1$.
Finally, again from the binomial identities,
\[
(1-t)Q^2(t)= \expec_t \biggl[\frac{m-S_m}{m} \bigl\{q
\bigl(Y_2,(m-1)/2\bigr)p\bigl(S_{m}-Y_2,
(m-1)/2\bigr)\mid S_{m} \bigr\} \biggr],
\]
with $\mathcal{L}(Y_2\mid S_{m}=k)=\Hypergeo(k,(m-1)/2,m-1)$,
$k=0,\ldots,m-1$.
\end{pf*}
\end{appendix}

\section*{Acknowledgements}
We are particularly grateful to the two referees and especially to the
Associate Editor, the comments and suggestions
have led to a very much improved version of this paper. We wish to thank
Libasse Coly for bringing the problem to our attention and
Johan Segers (Universit\'{e} Catholique de Louvain) for great advice
and interesting discussions. We also acknowledge
the financial support of the Natural Sciences and Engineering Research
Council of
Canada (NSERC).


%
%

\printhistory
\end{document}